\documentclass[a4paper,draft,12pt,reqno]{article}
\usepackage{enumerate}
\usepackage{amsmath,amssymb,amsthm,xspace,a4wide,amscd,latexsym,amscd}
\usepackage{xy}
\xyoption{all}
\DeclareMathOperator{\Hom}{Hom}
\DeclareMathOperator{\Ext}{Ext}
\newcommand{\bnabla}{\ensuremath{\overline{\nabla}}}
\newcommand{\bDelta}{\ensuremath{\overline{\Delta}}}

\DeclareMathOperator{\rad}{rad}

\newcommand{\add}{\operatorname{add}} 
 
\newcommand{\codim}{\operatorname{codim}} 
\newcommand{\tto}{\twoheadrightarrow} 

\renewcommand{\epsilon}{\varepsilon}

\renewcommand{\k}{\ensuremath{{\Bbbk}}}

\renewcommand{\phi}{\varphi}

\newcommand{\End}{\operatorname{End}}

\newcommand{\rank}{\operatorname{rank}}

\newcommand{\spar}{\operatorname{Tr}}
\newcommand{\DD}{\mathbb{D}}
\newcommand{\Der}{\mathcal{D}}
\newcommand{\Ka}{\mathcal{K}}

\begin{document} 
%\swapnumbers
\theoremstyle{plain}
%\numberwithin{subsection}{section}
%\newtheorem{thm}{Theorem}[section]
%\newtheorem{example}{Example}%[section]
%\newtheorem{definition}{Definition}%[section]
%\newtheorem{remark}{Remark}%[section]
\newtheorem{theorem}{Theorem}%[section]
\newtheorem{lemma}[theorem]{Lemma}%[section]
\newtheorem{proposition}[theorem]{Proposition}%[section]
\newtheorem{corollary}[theorem]{Corollary}%[section]

\title{Properly stratified algebras and tilting}
\author{Anders Frisk and Volodymyr Mazorchuk}
\date{}
\maketitle

\begin{abstract}
We study the properties of tilting modules in the context of properly
stratified algebras. In particular, we answer the question when the
Ringel dual of a properly stratified algebra is properly stratified
itself, and show that the class of properly stratified algebras for
which the characteristic tilting and cotilting modules coincide is 
closed under taking the Ringel dual. Studying stratified algebras, 
whose Ringel dual is properly stratified, we discover a new Ringel-type 
duality for such algebras, which we call the two-step duality. This 
duality arises from the existence of a new (generalized) tilting module 
for stratified algebras with properly stratified Ringel dual. We 
show that this new tilting module has a lot of interesting properties, 
for instance, its projective dimension equals the projectively defined 
finitistic dimension of the original algebra, it guarantees that the 
category of modules of finite projective dimension is contravariantly 
finite, and, finally, it allows one to compute the finitistic dimension of 
the original algebra in terms of the projective dimension of the 
characteristic tilting module.   
\end{abstract}

\section{Introduction}\label{s1} 

The concept of Ringel duality for quasi-hereditary algebras,
introduced in \cite{Ri}, has been shown to be a very important and 
useful tool for the study of different situations in which quasi-hereditary
algebras arise in a natural way, for example, for the study of the category
$\mathcal{O}$, Schur algebras and algebraic groups. 
The list of natural generalizations of quasi-hereditary algebras
starts with the so-called standardly stratified algebras, introduced
in \cite{CPS}, and the so-called properly stratified algebras,
introduced in \cite{Dl}. The study of the Ringel duality for
stratified algebras was originated in \cite{AHLU}, where the results,
analogous to those of Ringel, were obtained. Alternative approach to
the study of tilting modules for stratified algebras was developed in
\cite{PR,Xi}. 

The class of properly stratified algebras is a proper subclass of the
class of standardly stratified algebras. In particular, the results of
\cite{AHLU} perfectly apply to properly stratified algebras. However,
the very definition of the properly stratified algebras suggests that
such algebras must posses much more symmetric properties, for
example, because of the left-right symmetry of the class of such algebras. 
The principle motivation for the present paper was the question whether one
can obtain any nice new properties for the Ringel duals of properly stratified
algebras. 

Our first (relatively unexpected) observation is that the Ringel dual
of a properly stratified algebra does not have to be properly stratified, 
see the example in Subsection~\ref{s8.2}. This leads to the following natural 
question: when the Ringel dual of a stratified algebra is properly 
stratified? We answer this question in Section~\ref{s3} in terms of
the existence of special filtrations on tilting modules. In
Section~\ref{s4} we show that a natural class  of properly stratified
algebras is closed with respect to taking of the Ringel dual. This is 
the class of properly stratified algebras for which the characteristic 
tilting and cotilting modules coincide. This class appears naturally in 
\cite{MO} during the study of the finitistic dimension of stratified 
algebras and contains, in particular, all quasi-hereditary algebras. 

Assume now that $A$ is a stratified algebra, whose Ringel dual $R$ is
properly stratified. As a properly stratified algebra, the algebra $R$
has both tilting and cotilting modules. The classical Ringel duality
identifies the characteristic cotilting $R$-module with the
injective cogenerator of $A\mathrm{-mod}$, whereas the characteristic 
tilting $R$-module is identified with a possibly different $A$-module, 
which we call $H$. It happens that this module $H$ carries a lot of very 
important information about the algebra $A$ and the major part of our paper 
is devoted to the study of the properties of $H$. In particular,
in Sections~\ref{s5} and \ref{s6} we obtain the following:
\begin{enumerate}[(i)]
\item the module $H$ is a (generalized) tilting module;
\item the module $H$ is relatively injective with respect to the
subcategory of $A$-modules of finite projective dimension;
\item the projective dimension of $H$ equals the finitistic dimension of $A$;
\item existence of $H$ guarantees that the category of $A$-modules of
finite projective dimension is contravariantly finite in $A\mathrm{-mod}$;
\item if the algebra $A$ has a simple preserving duality, then the
existence of $H$ guarantees that the finitistic dimension of $A$ equals 
twice the projective dimension of the characteristic tilting $R$-module.
\end{enumerate}    
 
The classical tilting theory (see for example \cite{Ha}) motivates the
study of the endomorphism algebra of a (generalized) tilting module,
whenever one has such a module. Hence in Section~\ref{s7} we consider 
the opposite of the endomorphism algebra $B$ of the module $H$. We show 
that this algebra is standardly stratified and that its Ringel dual is properly 
stratified. This allows us to consider the corresponding module $H^{(B)}$. 
It happens that $H^{(B)}$ can be naturally identified with the injective 
cogenerator for $A^\mathrm{opp}\mathrm{-mod}$, which leads to a new 
covariant Ringel-type duality for stratified algebras having properly 
stratified Ringel dual. 

We show that this duality induces an equivalence between the category of 
all $A$-modules of finite projective dimension and the category of all
$B$-modules of finite injective dimension. It is not surprising that
this duality was not discovered in the theory of quasi-hereditary
algebras for in that case, and, more generally, in the case of properly 
stratified algebras for which the characteristic tilting and cotilting 
modules coincide, it degenerates to the identity functor. In particular, in
such case we always have $B=A$.

We finish the paper with various examples illustrating the new
duality, which include, in particular, certain categories of
Harish-Chandra bimodules over complex semi-simple Lie algebras, and
tensor products of quasi-hereditary and local algebras.

\section{Notation}\label{s2}

\subsection{General conventions}\label{s2.1} 

Let $A$ be a finite dimensional associative algebra with identity
over a field $\k$. For simplicity we assume that $\k$ is
algebraically closed. However, slightly modifying the dimension
arguments, one easily extends all the results to the case of an arbitrary
field. Denote by $\Lambda$ an index set for the isomorphism classes 
of simple $A$-modules, which we denote by $L(\lambda)$, $\lambda \in 
\Lambda$. We write $P(\lambda)$ for the projective cover
and $I(\lambda)$ for the injective hull of $L(\lambda)$. 

Throughout the paper a module means a {\em left} module. Our
main object will be some associative algebra $A$ (satisfying certain 
conditions, which we will specify later). Since we will consider more
than one algebra, to avoid confusion we adopt the following convention: 
$A$-modules will be written without any additional notation, for 
example $M$; modules over any other algebra, $B$ say, will be written as 
$M^{(B)}$ or ${}_B M$. By $A\mathrm{-mod}$ we denote the category of all 
finite-dimensional  left $A$-modules.      

Let $\mathcal{C}$ be a subclass of objects from $A\mathrm{-mod}$. Define 
$\mathcal{F}(\mathcal{C})$ as the full subcategory of $A\mathrm{-mod}$, which 
consists of all modules $M$ having a filtration, whose subquotients are isomorphic 
to modules from $\mathcal{C}$. Given an $A$-module, $M$, we define $\add(M)$ to 
be the full subcategory of $A\mathrm{-mod}$ containing all modules, which are 
isomorphic to direct summands of $M^k$ for some $k \geq 0$. Let $M$ and $N$ be two 
$A$-modules. We define the \textit{trace} $\spar_M(N)$ of $M$ in $N$ 
as the sum of images of all $A$-homomorphisms from $M$ to $N$.  

We denote by $\DD$ the usual duality functor
$\DD({}_-)=\Hom_{\k}({}_-,\k)$. 
Recall also that an algebra, $A$, has a \textit{simple preserving
duality} if there exists an exact contravariant and involutive equivalence 
${}^\circ:A\text{-mod} \to A\text{-mod}$, which preserves the isomorphism 
classes of simple $A$-modules.

For an algebra, $A$, we denote by $\mathcal{P}(A)^{<\infty}$ the full 
subcategory of $A$-mod, consisting of all $A$-modules of finite
projective dimension; and by $\mathcal{I}(A)^{<\infty}$ the full 
subcategory of $A$-mod, consisting of all $A$-modules of finite
injective dimension. We also denote by $\Der(A)$ the derived category
of $A\mathrm{-mod}$, by $\Der^b(A)$ its full subcategory consisting of all bounded 
complexes, and by $\Der^-(A)$ its full subcategory consisting of all right bounded 
complexes. By $\Ka(A\text{-mod})$ we denote the homotopy category of the 
category of all complexes of $A$-modules, by $\Ka^b(A\text{-mod})$ its full
subcategory consisting of all bounded complexes, and by $\Ka^-(A\text{-mod})$ 
its full subcategory consisting of all right bounded complexes.

Let $\mathcal{X}^{\bullet}$ be a complex in $\Ka(A\text{-mod})$ 
and $ j \in \mathbb{Z}$. We define the \textit{truncated complex} 
$t_{j}\mathcal{X}^{\bullet}$ to be the complex 
\begin{displaymath}
t_{j}\mathcal{X}^{\bullet}:\quad\quad \dots \rightarrow X_{j-2}
\overset{d_{j-2}}{\longrightarrow} X_{j-1}
\overset{d_{j-1}}{\longrightarrow} X_{j} 
\rightarrow 0 \rightarrow 0 \dots,
\end{displaymath}
where we keep the differentials $d_{i}$, $i<j$. A complex, 
$\mathcal{X}^{\bullet}=\{X_i:i\in\mathbb{Z}\}$, is called
\textit{positive} or \textit{negative} provided that $X_i=0$ for
all $i<0$ or $i>0$ respectively. For an $A$-module, $M$, we denote by 
$M^{\bullet}$ the corresponding complex in $\Ka^b(A\text{-mod})$, concentrated 
in degree zero.

\subsection{Stratified algebras and quasi-hereditary
algebras}\label{s2.2}  

Let $\leq$ be a linear order on $\Lambda$. For
$\lambda,\mu\in\Lambda$ we will write $\lambda<\mu$ provided 
that $\lambda\leq \mu$ and $\lambda\neq \mu$. For $\lambda \in 
\Lambda$ set $P^{>\lambda} = \oplus_{\mu>\lambda}P(\mu)$ and 
$I^{>\lambda} = \oplus_{\mu>\lambda}I(\mu)$. For each 
$\lambda \in \Lambda$ we define the \textit{standard module} 
\begin{displaymath}
\Delta(\lambda) = P(\lambda)/\spar_{P^{>\lambda}}(P(\lambda)),
\end{displaymath}
and the \textit{proper standard module} 
\begin{displaymath}
\bDelta(\lambda) =
\Delta(\lambda)/_{\sum_{f}\text{Im}f},  
\end{displaymath}
where the sum is taken over all
$f\in\rad\End_{A}\big(\Delta(\lambda)\big)$. Dually, we define the 
\textit{costandard module} 
\begin{displaymath}
\nabla(\lambda) = \bigcap_{f:I(\lambda) \to I^{>\lambda}}\text{Ker}f,
\end{displaymath}
and the \textit{proper costandard module} 
\begin{displaymath}
\bnabla(\lambda) = \bigcap_{f}
\text{Ker} f, 
\end{displaymath}
where the intersection is taken over all homomorphisms 
$f\in\rad\End_{A}\big(\nabla(\lambda)\big)$.

We can now define different classes of algebras which we will consider
in this paper. The pair $(A,\leq)$ is called a \textit{strongly 
standardly stratified algebra} or simply an \textit{SSS-algebra} if  
\begin{itemize}
\item[(SS)] the kernel of the canonical epimorphism $P(\lambda) 
\twoheadrightarrow \Delta(\lambda)$ has a filtration, whose  
subquotients are $\Delta(\mu)$ with $\lambda < \mu$. 
\end{itemize}
In several papers, see for example \cite{AHLU,ADL2}, the authors used the
name {\em standardly stratified algebras} for the algebras defined above.
This gives rise to a confusion with a more general original definition from
\cite{CPS}, where standardly stratified algebras were defined with respect
to a {\em partial pre-order} on $\Lambda$. Since in the present paper we 
will work with linear orders, we decided to use a different name.

The SSS-algebra $(A,\leq)$ is said to be \textit{properly stratified}
(see \cite{Dl}) if the following condition is satisfied:
\begin{itemize}
\item[(PS)] for each $\lambda \in \Lambda$ the module
$\Delta(\lambda)$ has a filtration, whose subquotients are isomorphic to 
$\bDelta(\lambda)$. 
\end{itemize}
Since the order $\leq$ will be fixed throughout the paper, we will usually 
omit it in the notation. It is easy to see (consult \cite{Dl}) that an 
SSS-algebra, $A$, is properly stratified if and only if $A^{\mathrm{opp}}$ is 
an SSS-algebra as well (with respect to the same order $\leq$). In particular, 
the algebra $A$ is properly stratified if and only if $A^{\mathrm{opp}}$ is.

The smallest class of algebras we will treat is the class of 
quasi-hereditary algebras, defined in the following way: 
the SSS-algebra $(A,\leq)$ is called \textit{quasi-hereditary} 
(see \cite{CPS0}) if the following condition is satisfied:
\begin{itemize}
\item[(QH)]  for each $\lambda \in \Lambda$ we have 
$\Delta(\lambda)=\bDelta(\lambda)$.  
\end{itemize}

If $A$ is standardly (properly) stratified then we denote by
$\mathcal{F}(\Delta)$ the category $\mathcal{F}(\mathcal{C})$, where 
$\mathcal{C}=\{\Delta(\lambda) \vert \lambda \in \Lambda \}$, and 
define $\mathcal{F}(\bDelta)$, $\mathcal{F}(\nabla)$ and
$\mathcal{F}(\bnabla)$ similarly.  

Let $A$ be an SSS-algebra. It was shown in \cite{AHLU} that the
category $\mathcal{F}(\Delta)\cap\mathcal{F}(\bnabla)$ is closed under 
taking direct summands, and that the indecomposable modules in this
category are indexed by $\lambda \in \Lambda$ in a natural way. 
The objects of $\mathcal{F}(\Delta)\cap\mathcal{F}(\bnabla)$ are
called \textit{tilting modules}. For $\lambda \in \Lambda$ we denote
by $T(\lambda)$ the (unique up to isomorphism) indecomposable object in 
$\mathcal{F}(\Delta)\cap\mathcal{F}(\bnabla)$ for which there exists
an exact sequence,
\begin{displaymath}
0\to \Delta(\lambda)\to T(\lambda)\to \mathrm{Coker}\to 0,
\end{displaymath}
such that $\mathrm{Coker}\in\mathcal{F}(\Delta)$ (see 
\cite[Lemma~2.5]{AHLU}). For $T = \oplus_{\lambda \in \Lambda}T(\lambda)$ 
we have $\mathcal{F}(\Delta)\cap\mathcal{F}(\bnabla) = \add(T)$. The module 
$T$ is usually called the \textit{characteristic tilting module}.

It also follows from \cite{AHLU} that, in the case when $A$ is 
properly stratified, the category 
$\mathcal{F}(\bDelta)\cap\mathcal{F}(\nabla)$ is closed under taking 
direct summands, and that the indecomposable modules in this
category are indexed by $\lambda \in \Lambda$ in a natural way. 
The objects of $\mathcal{F}(\bDelta)\cap\mathcal{F}(\nabla)$ are
called \textit{cotilting modules}.  For $\lambda \in \Lambda$ we 
denote by $C(\lambda)$ the (unique up to isomorphism) indecomposable object in 
$\mathcal{F}(\bDelta)\cap\mathcal{F}(\nabla)$ for which there exists an exact
sequence,
\begin{displaymath}
0\to \mathrm{Ker}\to C(\lambda)\to \nabla(\lambda)\to 0,
\end{displaymath}
such that $\mathrm{Ker}\in\mathcal{F}(\nabla)$. For 
$C=\oplus_{\lambda \in \Lambda}C(\lambda)$ we have 
$\mathcal{F}(\bDelta)\cap\mathcal{F}(\nabla) = \add(C)$. The module 
$C$ is usually called the \textit{characteristic cotilting module}. 

We would like to point out that the name {\em (co)tilting module} introduced above
is slightly confusing. In the terminology of the classical tilting theory
the (co)tilting modules as defined above are only {\em partial (co)tilting modules}. 
On the other hand, $A\mathrm{-mod}$ usually contains many other modules,
which are (co)tilting in the classical sense but are not related to the stratified
structure and thus to the (co)tilting modules defined above.
However, we will keep the above name since it is now standard and commonly
accepted in the theory of quasi-hereditary and stratified algebras. Later on,
in Subsection~\ref{s6.1} we will recall the definition of a (generalized)
tilting modules from the classical tilting theory. A tilting module as defined
above is a (generalized) tilting module in the classical sense if and only if
it contains a direct summand, which is isomorphic to the characteristic tilting 
module.

We set 
$L = \oplus_{\lambda\in \Lambda}L(\lambda)$, 
$\Delta = \oplus_{\lambda\in \Lambda}\Delta(\lambda)$, 
$\nabla = \oplus_{\lambda\in \Lambda}\nabla(\lambda)$, 
$\bDelta = \oplus_{\lambda\in \Lambda}\bDelta(\lambda)$,  
$\bnabla = \oplus_{\lambda\in \Lambda}\bnabla(\lambda)$,   
$I = \oplus_{\lambda\in \Lambda}I(\lambda)$, and  
$P = \oplus_{\lambda\in \Lambda}P(\lambda)$. 

Throughout the paper we will multiply the maps from the left to the
right, for instance the composition $g\circ f$ of the map $f:X\to Y$ 
and the map $g:Y\to Z$ is denoted by $fg$. 

\subsection{Ringel dual for standardly stratified algebras}\label{s2.3} 

Let $(A,\leq)$ be an SSS-algebra. Following \cite{AHLU} we define the 
\textit{Ringel dual} $R=R(A)$ of $A$ as the algebra 
$R=\text{End}_{A}(T)$. This algebra comes together with the 
\textit{Ringel duality functor}
\begin{equation}\label{funktor}
F=F^{A}:A\text{-mod} \rightarrow R\text{-mod,} \quad
\text{ defined by } \quad F({}_-) = \text{Hom}_{A}(T,{}_-).
\end{equation}

Due to \cite[Theorem 2.6]{AHLU} the algebra $R^{opp}$ is an
SSS-algebra with the same indexing set $\Lambda$, but with respect to
the order $\leq_{R}$, which is opposite to the original order $\leq$. 
The Ringel duality asserts that the
algebra $A^{opp}$ is Morita equivalent to the Ringel dual of the algebra
$R^{opp}$. Moreover, the functor $F$ sends $\bnabla$ to
$\bDelta^{(R)}$ (the latter being defined with respect to $\leq_{R}$) 
and induces an equivalence between $\mathcal{F}(\bnabla)$ and 
$\mathcal{F}(\bDelta^{(R)})$.  

\section{When is the Ringel dual properly stratified?}\label{s3} 

Let $(A,\leq)$ be an SSS-algebra. First of all we would like to
determine when the Ringel dual $R$ of $A$ is properly stratified. 

For $\lambda \in \Lambda$ set $T^{<\lambda} = 
\oplus_{\mu<\lambda}T(\mu)$. Define further
$S(\lambda)=\spar_{T^{<\lambda}}(T(\lambda))$ and put 
$N(\lambda)=T(\lambda)/S(\lambda)$ giving the following short exact sequence: 
\begin{equation}\label{seq1}
0 \rightarrow S(\lambda) \rightarrow T(\lambda) \rightarrow
N(\lambda) \rightarrow 0.
\end{equation}
Finally, set $N =\oplus_{\lambda \in \Lambda}N(\lambda)$ and $\mathcal{F}(N)
=\mathcal{F}(\mathcal{C})$, where $\mathcal{C}=\{N(\lambda) \vert
\lambda \in \Lambda \}$. We start with another description of $S(\lambda)$.

\begin{lemma}\label{propSL}
For every $\lambda \in \Lambda$ the module $S(\lambda)$ is the unique
submodule $M$ of $T(\lambda)$ which is characterized by the following
properties:
\begin{enumerate}[(a)]
\item\label{propSL1} $M\in\mathcal{F}(\{\bnabla(\mu) | \mu <
\lambda\})$.   
\item\label{propSL2}
$T(\lambda)/M\in\mathcal{F}(\{\bnabla(\lambda)\})$.  
\end{enumerate}
\end{lemma}

\begin{proof}
Let 
\begin{displaymath}
0=T_0\subset T_1 \subset \dots \subset T_k=T(\lambda)
\end{displaymath}
be a proper costandard filtration of $T(\lambda)$. Using 
\begin{displaymath}
\Ext_A^1\big(\bnabla(\mu_1),\bnabla(\mu_2)\big) = 0 
\quad\text{ for all }\mu_1<\mu_2,
\end{displaymath}
we can assume that there exists $0 \leq l < k$ such that $T_i/T_{i-1} 
\cong \bnabla(\mu)$ with $\mu<\lambda$ for all $i\leq l$,
and $T_i/T_{i-1} \cong \bnabla(\lambda)$  for all $i> l$. For $\mu <
\lambda$ we apply the functor $\Hom_A(T(\mu),{}_-)$ to 
the short exact sequence
\begin{displaymath}
0 \to T_l \to T(\lambda) \to \text{Coker} \to 0,
\end{displaymath}
and obtain 
\begin{displaymath}
0 \to\Hom_A(T(\mu),T_l) \to\Hom_A(T(\mu),T(\lambda)) \to 
\Hom_A(T(\mu),\text{Coker}).
\end{displaymath}
As $\text{Coker}\in\mathcal{F}(\{\bnabla(\lambda)\})$, 
$T(\mu)\in\mathcal{F}(\Delta)$ and $[T(\mu):\Delta(\lambda)]=0$ we
get $\Hom_A(T(\mu),\text{Coker})=0$. This implies that 
$\Hom_A(T(\mu),T_l) = \Hom_A(T(\mu),T(\lambda))$
and thus the image of any map $f:T(\mu)\to T(\lambda)$ belongs to
$T_l$. In particular, $S(\lambda)\subset T_l$. 

To prove that $S(\lambda)= T_l$ we show by induction on $i$ that for
every $0\leq i\leq l$ the module $T_i$ is a quotient of some
$T^{(i)}\in \add\left(\oplus_{\nu<\lambda}T(\nu)\right)$. The statement
is obvious for $i=0$, and thus we have to prove only the induction step 
$i\implies i+1$. Consider the short exact sequence 
\begin{displaymath}
0 \to T_i \to T_{i+1} \to \bnabla(\nu) \to 0,
\end{displaymath}
where $\nu<\lambda$. Then we have the map $T^{(i)}\tto T_i \to 
T_{i+1}$, and an epimorphism, $T(\nu)\tto \bnabla(\nu)$. Since 
$T_i, T_{i+1}, \bnabla(\nu)\in\mathcal{F}(\bnabla)$, and
$T(\nu)\in\mathcal{F}(\Delta)$, the epimorphism
$T(\nu)\tto \bnabla(\nu)$ can be lifted to a map,  $T(\nu)\to T_{i+1}$,
giving an epimorphism, $T^{(i)}\oplus T(\nu)\tto T_{i+1}$. Therefore we
finally obtain $S(\lambda)= T_l$. In particular, we see that
$S(\lambda)$ satisfies both \eqref{propSL1} and \eqref{propSL2},
moreover, it is unique by construction. This completes the proof.
\end{proof}

From the proof of Lemma~\ref{propSL} we also obtain the following
information: 

\begin{corollary}\label{propertiesofN}
\begin{enumerate}[(1)]
\item\label{hlop1} For each $\lambda \in \Lambda$ the module $N(\lambda)$ has a
filtration, whose subquotients are isomorphic to $\bnabla(\lambda)$.
\item\label{hlop2} If $\mu < \lambda$ then $\text{Hom}_{A}(T(\mu),N(\lambda)) = 0$.
\end{enumerate}
\end{corollary}

Now we are ready to formulate and prove a characterization of those 
SSS-algebras, whose Ringel dual is properly stratified. 

\begin{theorem}\label{Ringelpropstr}
Let $(A,\leq)$ be an SSS-algebra. Then the following
assertions are equivalent:
\begin{enumerate}[(I)] 
\item\label{Rin1} The Ringel dual $(R,\leq_{R})$ is properly
stratified. 
\item\label{Rin2} For each $\lambda \in \Lambda$ we have $S(\lambda)
\in \mathcal{F}(N)$.
\item\label{Rin3} For each $\lambda \in \Lambda$ we have $T(\lambda)
\in \mathcal{F}(N)$.
\end{enumerate}
\end{theorem}

\begin{proof} 
The equivalence of the conditions \eqref{Rin2} and \eqref{Rin3} follows 
from \eqref{seq1}.

Since the functor $F$ induces an exact equivalence from
$F(\bnabla)$ to $F(\bDelta^{(R)})$, the short exact
sequence \eqref{seq1} gives the short exact sequence 
\begin{equation}\label{nereiR}
0 \rightarrow F(S(\lambda)) \rightarrow F(T(\lambda)) \rightarrow  
F(N(\lambda)) \rightarrow 0.
\end{equation}     
Moreover, $F(T(\lambda)) = P^{(R)}(\lambda)$ by definition.

\eqref{Rin2}$\implies$\eqref{Rin1}. The algebra $R^{opp}$ is an
SSS-algebra by \cite{AHLU}, hence it is enough to prove that so is
$R$. We start with the following useful statement, which we will
also use later on.

\begin{lemma}\label{lref1}
Let $(A,\leq)$ be an SSS-algebra. Then $F(N(\lambda)) 
\cong\Delta^{(R)}(\lambda)$.
\end{lemma}

\begin{proof}
Using Lemma~\ref{propSL}, the exactness of $F$ on $\mathcal{F}(\bnabla)$,
and $F(\bnabla(\nu))=\bDelta^{(R)}(\nu)$ for all $\nu\in\Lambda$
we obtain that $F(S(\lambda))$ has a filtration with subquotients 
$\bDelta^{(R)}(\mu)$, $\mu < \lambda$. This implies that 
$F(S(\lambda))$ includes into the trace of $\oplus_{\nu>_{R}\lambda}F(T(\nu))$ in 
$F(T(\lambda))$. On the other hand, from Corollary~\ref{propertiesofN}\eqref{hlop2} 
it follows that the latter inclusion is in fact the equality. Now the statement
follows from the definition of a standard module.
\end{proof}

Exactness of $F$ on $\mathcal{F}(\bnabla)$ and Lemma~\ref{lref1}
imply that $P^{(R)}\in \mathcal{F}(\{F(N(\lambda))|\lambda \in \Lambda\})$.
This means that $R$ is an SSS-algebra, and hence is properly stratified,
completing the proof of the implication \eqref{Rin2}$\implies$\eqref{Rin1}.
 
\eqref{Rin1}$\implies$\eqref{Rin2}. Assume that $R$ is properly
stratified, in particular is an SSS-algebra. For every $\lambda \in
\Lambda$ consider the short exact sequence 
\begin{displaymath}
0 \to K^{(R)}(\lambda) \to P^{(R)}(\lambda) \to \Delta^{(R)}(\lambda)
\to 0.
\end{displaymath}
Conditions (SS) and (PS) ensure that $P^{(R)}(\lambda)\in 
\mathcal{F}(\bDelta^{(R)})$, $\Delta^{(R)}(\lambda)\in 
\mathcal{F}(\{\bDelta^{(R)}(\lambda)\})$, $K^{(R)}(\lambda) \in 
\mathcal{F}(\{\bDelta^{(R)}(\mu)|\mu>_R \lambda\})$. Now we can 
apply $F^{-1}$, which is exact on  $\mathcal{F}(\bDelta^{(R)})$,
and obtain the exact sequence 
\begin{displaymath}
0 \to F^{-1}\big(K^{(R)}(\lambda)\big) \to T(\lambda) 
\to F^{-1}\big(\Delta^{(R)}(\lambda)\big) \to 0,
\end{displaymath}
moreover, 
$F^{-1}\big(\Delta^{(R)}(\lambda)\big)\in\mathcal{F}(\{\bnabla(\lambda)\})$, 
and $ F^{-1}\big(K^{(R)}(\lambda)\big)\in \mathcal{F}(\{\bnabla(\mu)|\mu<
\lambda\})$. From Lemma~\ref{propSL} we obtain that   
$S(\lambda)=F^{-1}\big(K^{(R)}(\lambda)\big)$ and 
$N(\lambda)=F^{-1}\big(\Delta^{(R)}(\lambda)\big)$. 
Using $K^{(R)}(\lambda)\in\mathcal{F}(\Delta^{(R)})$ we get 
$S(\lambda)\in F^{-1}\left(\mathcal{F}(\Delta^{(R)})\right)=\mathcal{F}(N)$.
This proves the implication \eqref{Rin1}$\implies$\eqref{Rin2} and thus the
proof is complete.
\end{proof}

\section{The Ringel dual of a properly stratified algebras}\label{s4} 

In general the Ringel dual of a properly stratified algebra does not
need to be properly stratified, see the example in Subsection~\ref{s8.2}. In
this section we show that for one natural class of properly stratified
algebras this can always be guaranteed (in fact, we prove even more,
namely, that this class is closed under taking the Ringel dual). The
class consists of all properly stratified algebras for which the 
characteristic tilting and cotilting modules coincide. This class appeared 
in \cite{MO} where it was shown that the finitistic dimension of an 
algebra with duality from this class is twice the projective dimension of the 
characteristic tilting module (we will extend this result in Subsection~\ref{s6.4}). 
This class contains, in particular, all quasi-hereditary algebras.       

\begin{theorem}\label{thmringel}
Assume that $(A,\leq)$ is a properly stratified algebra for which the 
characteristic tilting and cotilting modules coincide. Then the Ringel dual
$(R,\leq_{R})$ is properly stratified, moreover, the characteristic tilting 
and cotilting $R$-modules coincide as well.  
\end{theorem}

\begin{proof}
Let $(A,\leq)$ be a properly stratified algebra. By
Theorem~\ref{Ringelpropstr}, to show that $(R,\leq_{R})$ is
properly stratified it is enough to show that
$S(\lambda) \in \mathcal{F}(N)$ for each $\lambda \in
\Lambda$.  Recall that $T(\lambda) = C(\lambda)$ by
our assumption. This gives us the following short exact sequence:
\begin{displaymath}
0\to\text{Ker}\to T(\lambda)\to\nabla(\lambda)\to 0, 
\end{displaymath}
where $\text{Ker}$ has a filtration by $\nabla(\mu)$, $\mu <
\lambda$. Since $A$ is properly stratified, filtrations by costandard
modules extend to filtrations by the corresponding proper costandard
modules. Hence we can use Lemma~\ref{propSL} to obtain
$S(\lambda) =\text{Ker}$ and $N(\lambda)
=\nabla(\lambda)$. This implies that $(R,\leq_{R})$ is
properly stratified.

It is left to show that $T^{(R)} = C^{(R)}$. By
\cite[Theorem~2.6(vi)]{AHLU}, $F(I)$ is the characteristic
cotilting $R$-module, and, since $R$ is properly stratified,
we need only to show that $F(I)\in\mathcal{F}(\Delta^{(R)})$.
However, since $A$ is properly stratified, we have
$I\in\mathcal{F}(\nabla)=\mathcal{F}(N)$, $F$ is exact on
$\mathcal{F}(\nabla)$, and $F(N)=\Delta^{(R)}$ by Lemma~\ref{lref1}.
This implies that $F(I)\in\mathcal{F}(\Delta^{(R)})$ and completes the proof.
\end{proof}

\section{Module $H$: definition and basic properties}\label{s5} 

Now we can make the {\bf principal assumption} until the end of the paper:
{\em we assume that $(A,\leq)$ is an SSS-algebra such that the Ringel dual 
$(R,\leq_{R})$ is properly stratified}. 

There are two possibilities. The first one is the case when $T^{(R)}= C^{(R)}$. 
In this case Theorem~\ref{thmringel} says that the algebra $A$ itself is properly
stratified and, moreover, that $T= C$. Hence the behavior of the
Ringel duality in this case is completely similar to the classical
case of quasi-hereditary algebras (which is in fact a special case of
this situation). 

The second case is when $T^{(R)} \not\cong C^{(R)}$. In this case it 
is possible that the algebras $\End_R(T^{(R)})$ and $\End_R(C^{(R)})$ 
are quite different (see the example in Subsection~\ref{s8.2}). Later on
in the paper we will show that this situation leads to a new Ringel  
type duality, which we call the \textit{two-step duality}, 
on a subclass of SSS-algebras. In the special case of quasi-hereditary 
algebras, and even in the more general situation of properly stratified algebras 
for which the characteristic tilting and cotilting modules coincide, the two-step
duality degenerates to the identity functor.  

Since the Ringel duality functor $F$ induces an equivalence between
$\mathcal{F}(\bnabla)$ and $\mathcal{F}(\bDelta^{(R)})$, we can
take the module $T^{(R)}(\lambda)\in\mathcal{F}(\bDelta^{(R)})$ and
define $H(\lambda) = F^{-1}\big(T^{(R)}(\lambda)\big)$ for all 
$\lambda\in \Lambda$. Set $H = \oplus_{\lambda \in
\Lambda}H(\lambda)$. The module $H$ will be the main object of
our  interest in the sequel. We start with some basic properties of 
the modules $N$ and $H$.

\begin{proposition}\label{twostep}
\begin{enumerate}[(i)] 
\item\label{H1} All tilting modules belong to $\mathcal{F}(N)$ and are
exactly the relatively (ext-)projective modules in $\mathcal{F}(N)$, that is
for $M \in \mathcal{F}(N)$ the module $M$ is tilting if and only if 
$\text{Ext}_{A}^{i}(M,N)=0$ for all  $i>0$.
\item\label{H2} $\add(H)\subset \mathcal{F}(N)$ and the modules from
$\add(H)$ are exactly the relatively (ext-)in\-jec\-tive modules in 
$\mathcal{F}(N)$, that is for $M \in \mathcal{F}(N)$ we have
$M\in \add(H)$ if and only if $\text{Ext}_{A}^{i}(N,M)=0$ for all $i>0$.
\end{enumerate}
\end{proposition}

\begin{proof}
We start with \eqref{H1}. That all tilting modules belong to 
$\mathcal{F}(N)$ follows from Theorem~\ref{Ringelpropstr}\eqref{Rin3}.
Assume that $M$ is a tilting module. Then 
$\Ext_A^i(M,N)=0$ for all $i>0$ by \cite[Lemma~1.2]{AHLU} 
since $M\in\mathcal{F}(\Delta)$ and
$N\in\mathcal{F}(\bnabla)$. To prove the opposite
statement, assume that $M\in \mathcal{F}(N)$ and 
$\text{Ext}_{A}^{i}(M,N)=0$, for all
$i>0$. Since $M\in \mathcal{F}(N) \subset
\mathcal{F}(\bnabla)$, we only need to prove that 
$M\in \mathcal{F}(\Delta)$. 

Let $\lambda \in \Lambda$ and $X$ be a submodule of 
$N(\lambda)$ such that we have the following short exact 
sequence:
\begin{displaymath} 
0 \to X \to N(\lambda) \to \bnabla(\lambda)  
\to \ 0, 
\end{displaymath}
where $X$ has a filtration with subquotients
$\bnabla(\lambda)$. Applying  
$\Hom_A(M,{}_-)$ to this short exact sequence we get the 
following fragment in the long exact sequence  
\begin{displaymath} 
\dots \to \Ext_A^i(M,N(\lambda))
\to \Ext_A^i(M,\bnabla(\lambda)) 
\to \Ext_A^{i+1}(M,X) \to
\Ext_A^{i+1}(M,N(\lambda)) \dots.
\end{displaymath} 
The condition $\Ext_A^i(M,N(\lambda)) = 0$, $i >0$, gives 
us a dimension shift between the spaces
$\Ext_A^i(M,\bnabla(\lambda))$ and
$\Ext_A^{i+1}(M,X)$. 

To proceed we need the following statement:

\begin{lemma}\label{projNN}
The module $N$ has finite projective dimension.
\end{lemma}

\begin{proof}
We prove that $\text{p.d.}(N(\lambda))<\infty$ by induction on
$\lambda \in \Lambda$. Suppose $\lambda$ is minimal. Then
$N(\lambda) \cong T(\lambda)$ and hence
$\text{p.d.}(N(\lambda))<\infty$. Now assume by induction
that for all  $\mu < \lambda$ we have
$\text{p.d.}(N(\mu))<\infty$. Since $S(\lambda)$ is filtered by
$N(\mu)$, $\mu < \lambda$, it follows that 
$\text{p.d.}(S(\lambda))<\infty$. The exact sequence \eqref{seq1}
and $\text{p.d.}(T(\lambda))<\infty$ now implies 
$\text{p.d.}(N(\lambda))<\infty$ and completes the proof.
\end{proof}

Since $M\in\mathcal{F}(N)$, Lemma~\ref{projNN} implies that
$\text{p.d.}(M)<\infty$. This forces the equality 
$\Ext_A^i(M,\bnabla(\lambda))=0$ for all $i$ big
enough. Thus, the dimension shift and the fact that $X$ is
filtered by $\bnabla(\lambda)$, guarantees that
$\Ext_A^i(M,\bnabla(\lambda))=0$ for all $i>0$. From
\cite[Theorem~1.6]{AHLU}  we now obtain  
$M\in \mathcal{F}(\Delta)$, which completes the proof of \eqref{H1}.

Now let us prove \eqref{H2}. 
Since $F:\mathcal{F}(\bnabla)\to\mathcal{F}(\bDelta^{(R)})$ is
an equivalence with $\mathcal{F}(\bnabla)=\{X\in A\mathrm{-mod}|
\text{Ext}_{A}^{>0}\left(T,M\right)=0\}$ and $\mathcal{F}(N)\subset
\mathcal{F}(\bnabla)$, we obtain that 
for all $i\geq  0$ and for all $M_1, M_2 \in
\mathcal{F}(N)$ we have 
\begin{equation}\label{niceformula} 
\text{Ext}_{A}^{i}\left(M_1,M_2\right) \cong
\text{Ext}_{R}^{i}\left(F\big(M_1\big),F\big(M_2\big)\right).
\end{equation}
Put $M_1 = N$ and $M_2 = H$. By Lemma~\ref{lref1} we 
have $F\big(N\big)=\Delta^{(R)}$ and by the definition
we have $F\big(H\big)=T^{(R)}$. Therefore 
$H\in \mathcal{F}(N)$ and \eqref{niceformula} guarantees
\begin{displaymath} 
\text{Ext}_{A}^{i}\left(N,H\right) \cong
\text{Ext}_{R}^{i}\left(\Delta^{(R)},T^{(R)}\right) = 0,
\end{displaymath}
as $T^{(R)}\in\mathcal{F}(\bnabla^{(R)})$. This proves that $H$
is relatively injective in $\mathcal{F}(N)$. 

To prove the opposite statement, we assume that
$M\in\mathcal{F}(N)$ and
$\text{Ext}_{A}^{i}(N,M)=0$ for all $i>0$. Using 
\eqref{niceformula} we get that
$\text{Ext}_{R}^{i}\left(\Delta^{(R)},F\big(M\big)\right) = 0$, 
which implies that
$F\big(M\big)\in\mathcal{F}(\bnabla^{(R)})$ by
\cite[Theorem~1.6]{AHLU}. Moreover, 
$F\big(M\big)\in \mathcal{F}(\Delta^{(R)})$ since
$M\in\mathcal{F}(N)$, $F(N)=\Delta^{(R)}$, and $F$ is exact on 
$\mathcal{F}(N)\subset\mathcal{F}(\bnabla)$. 
Hence $F\big(M\big)$ is a tilting $R$-module and thus $M\in \add(H)$
by the definition of $H$. This completes the proof.
\end{proof}

\begin{proposition}\label{prref3}
Let $\lambda\in\Lambda$. Then there exist the following exact sequences:
\begin{eqnarray}
0\to \bnabla(\lambda)\to H(\lambda)\to \mathrm{Coker}_1\to 0,
\label{prref3.1} \\
0\to N(\lambda)\to H(\lambda)\to \mathrm{Coker}_2\to 0,
\label{prref3.2} 
\end{eqnarray}
where $\mathrm{Coker}_1$ has a filtration with subquotients
$\bnabla(\mu)$, $\mu\geq \lambda$ and 
$\mathrm{Coker}_2$ has a filtration with subquotients
$N(\mu)$, $\mu> \lambda$.
\end{proposition}

\begin{proof}
Applying $F^{-1}$ to the exact sequence
\begin{equation}\label{prref3.3}
0\to \Delta^{(R)}(\lambda)\to T^{(R)}(\lambda)\to \widetilde{\mathrm{Coker}}^{(R)}_2\to 0,
\end{equation}
where $\widetilde{\mathrm{Coker}}^{(R)}_2\in\mathcal{F}(\Delta^{(R)})$ (see 
\cite[Lemma~2.5(iv)]{AHLU}), gives \eqref{prref3.2}. 

Since $R$ is properly stratified, all standard modules
have proper standard filtrations and hence \eqref{prref3.3} gives rise to
the exact sequence
\begin{equation}\label{prref3.4}
0\to \bDelta^{(R)}(\lambda)\to T^{(R)}(\lambda)\to \widetilde{\mathrm{Coker}}^{(R)}_1\to 0,
\end{equation}
where $\widetilde{\mathrm{Coker}}^{(R)}_1\in\mathcal{F}(\bDelta^{(R)})$. Applying
$F^{-1}$ to \eqref{prref3.4} gives \eqref{prref3.1}. This completes the proof.
\end{proof}

\section{Module $H$: advanced properties}\label{s6} 

\subsection{$H$ is a (generalized) tilting module}\label{s6.1} 

Recall that a module, $M$, over an associative algebra $A$ is called a  
\textit{(generalized) tilting module} (see for example 
\cite[Chapter~III]{Ha}) if the following three conditions are satisfied:    
\begin{enumerate}[(i)] 
\item\label{Mprop1} $\text{Ext}_{A}^{i}(M,M) = 0$, $i>0$;   
\item\label{Mprop2} $\text{p.d.}(M)<\infty$;   
\item\label{Mprop3} there is an exact sequence   
$0 \rightarrow {}_{A}A \rightarrow M_{0} \rightarrow M_{1} 
\rightarrow \dots \rightarrow M_{k} \rightarrow 0$, where $k\geq 0$
and $M_{i}\in \add(M)$.    
\end{enumerate}
We would like to emphasize once more that a tilting module, $M$, over an 
SSS-algebra is a (generalized) tilting module in the above sense if
and only if $M$ contains a direct summand, which is isomorphic to the
characteristic tilting module, see Subsection~\ref{s2.2}. 
Now we can state the following:

\begin{proposition}\label{Histilting}
The module $H$ is a (generalized) tilting module. 
\end{proposition}

\begin{proof}
We split the proof into a sequence of lemmas.

\begin{lemma}\label{extHH}
\begin{displaymath}
\Ext^{i}(H,H)=0 \quad\quad \text{ for all } i> 0.
\end{displaymath}
\end{lemma}

\begin{proof}
This follows directly from the second statement of 
Proposition~\ref{twostep}.
\end{proof}

\begin{lemma}\label{projN}
The module $H$ has finite projective dimension.
\end{lemma}

\begin{proof}
Since $H \in\mathcal{F}(N)$, the statement follows 
from  Lemma~\ref{projNN}.
\end{proof}

\begin{lemma}\label{Hres}
For every $X \in \add(T)$ there exists a minimal coresolution,
\begin{displaymath}
0 \rightarrow X \rightarrow H_{0} 
\rightarrow H_{1} \rightarrow \dots
\rightarrow H_{k} \rightarrow 0,
\end{displaymath}
where $H_{i} \in \add(H)$ for all $i$ and $ 0 \leq k 
\leq \text{p.d.}(T^{(R)})$.
\end{lemma}

Note that the length of the coresolution in Lemma~\ref{Hres} is
estimated in terms of the projective dimension of the tilting module 
over the Ringel dual $R$.

\begin{proof}
The module $F(X)$ is a projective $R$-module and hence, by
\cite[Section~III.2.2]{Ha}, there exists a finite coresolution,
\begin{equation}\label{eq2003}
0 \to F\big(X\big) \to T^{(R)}_0 \to T^{(R)}_1 \to \dots 
\to T^{(R)}_k \to 0,
\end{equation}
where $T^{(R)}_i \in \add(T^{(R)})$ for all $i$, of length
$k \leq \text {p.d.}(T^{(R)})$. Applying $F^{-1}$ to \eqref{eq2003} 
yields the exact sequence
\begin{displaymath}  
0 \to X \to F^{-1}(T^{(R)}_0) \to F^{-1}(T^{(R)}_1) \to \dots \to 
F^{-1}(T^{(R)}_k) \to 0, 
\end{displaymath}
where $F^{-1}(T^{(R)}_i) \in \add(H)$ for all $i$.  
This completes the proof.
\end{proof}

\begin{lemma}\label{tiltingresofA}
There is an exact sequence 
\begin{displaymath}
0 \rightarrow {}_{A}A \rightarrow H_{0} \rightarrow H_{1} \rightarrow
\dots \rightarrow H_{m} \rightarrow 0,
\end{displaymath}
where $H_{i} \in \add(H)$. 
\end{lemma}

\begin{proof}
Choose a minimal tilting coresolution, 
\begin{displaymath}
0 \rightarrow {}_A A \rightarrow T_{0} \rightarrow T_{1} \rightarrow
\dots \rightarrow T_k \rightarrow 0,
\end{displaymath}
where $T_{i} \in \add(T)$ and $ k = \text{p.d.}(T)$, and consider the
corresponding positive complex 
\begin{equation}\label{FF1}
\dots \to 0 \to T_{0} \rightarrow T_{1} \rightarrow \dots
\rightarrow T_k \rightarrow 0 \to \dots,
\end{equation}
in $\Ka^b(\add(T))$, 
whose only non-zero homology is in degree zero and equals ${}_A A$. 
Applying Lemma~\ref{apl1} to $\mathtt{A}=A$, $X^{(\mathtt{A})}=T$, 
$Y^{(\mathtt{A})}=H$ and the complex \eqref{FF1} we obtain a complex,
\begin{displaymath}
\dots \to 0 \to H_{0} \rightarrow H_{1} \rightarrow \dots
\rightarrow H_{k'} \rightarrow 0 \to \dots,
\end{displaymath}
in $\Ka^b(\add(H))$, which is quasi-isomorphic to \eqref{FF1}. This 
completes the proof.
\end{proof}

The proof of Proposition~\ref{Histilting} is now also complete.
\end{proof}

We remark that, since $H$ is a (generalized) tilting module, the
minimal length of the coresolution, given by 
Lemma~\ref{tiltingresofA}, is equal to $\text{p.d.}(H)$, see for
example \cite[Chapter~III]{Ha}.

\subsection{$H$ and $\text{fin.dim.}(A)$}\label{s6.2} 

Recall that the (projectively defined) finitistic dimension of $A$ is
defined as follows:
\begin{displaymath}
\text{fin.dim}(A)=\text{sup}\{\text{p.d.}(M)
\vert M \in A\text{-mod},\,\, \text{p.d.}(M)<\infty\}.
\end{displaymath} 

In \cite{AHLU2} it was shown that the projectively defined finitistic 
dimension of a properly stratified (and even SSS-) algebra is 
finite (and even that the injectively defined finitistic dimension of 
such algebra is finite as well). In this subsection we show that the module 
$H$ can be used to effectively compute $\text{fin.dim}(A)$.

\begin{proposition}\label{coresofM}
Let $M \in A\mathrm{-mod}$ and $\text{p.d.}(M)<\infty$. Then
there exists a finite coresolution, 
\begin{displaymath}
0 \to M \to H_0 \to \dots \to H_k \to 0,
\end{displaymath}
where $H_i \in \add(H)$ and $k\geq 0$.
\end{proposition}

\begin{proof}
First we choose a minimal projective resolution,
\begin{displaymath}
0 \to P_n \to \dots \to P_1 \to P_0 \to M \to 0,
\end{displaymath}
of $M$, and obtain the complex 
\begin{displaymath}
\mathcal{P}^{\bullet}:\quad\quad \dots \to 0 \to P_n \to 
\dots \to P_1 \to P_0 \to 0 \to \dots,
\end{displaymath}
where $P_0$ stays in degree zero. The complex $\mathcal{P}^{\bullet}$ 
has the only non-zero homology in degree zero, which equals $M$.
Applying Lemma~\ref{apl1} to $\mathtt{A}=A$, $X^{(\mathtt{A})}=P$, $Y^{(\mathtt{A})}=H$ 
and  the complex $\mathcal{P}^{\bullet}[-n]$ we obtain a complex, shifting
which by $n$ in the derived category gives the complex 
\begin{displaymath}
\mathcal{H}^{\bullet}:\quad\quad \dots \to 0 \to H_{-n} 
\to \dots \to H_{-1} \to H_0 \to H_1 \to \dots
\to H_k \to 0 \to \dots,
\end{displaymath}
in $\Ka^b(\add(H))$, where $k\geq 0$, which is quasi-isomorphic
to $\mathcal{P}^{\bullet}$. This implies that the only non-zero homology of
$\mathcal{H}^{\bullet}$ is in degree zero and equals $M$.
Let us show that $\mathcal{H}^{\bullet}$ is quasi-isomorphic to a
positive complex from $\Ka^b(\add(H))$. In fact we will show by a
downward induction on $l$ that for every $0 \leq l \leq n$ there exists a
complex,  
\begin{displaymath}
\mathcal{H}(l)^{\bullet}: \dots \to 0 \to H(l)_{-l} 
\to \dots \to H(l)_{-1} \to H(l)_0 \to H(l)_1 \to \dots
\to H(l)_k \to 0 \to \dots,
\end{displaymath}
in $\Ka^b(\add(H))$, which is quasi-isomorphic to
$\mathcal{H}^{\bullet}$.

For $l=n$ the statement is obvious with 
$\mathcal{H}(n)^{\bullet}=\mathcal{H}^{\bullet}$.
Assume that the complex $\mathcal{H}(l)^{\bullet}$,  $0<l\leq n$, is 
constructed and let us construct the complex
$\mathcal{H}(l-1)^{\bullet}$. 

From $\mathcal{H}(l)^{\bullet}$ we obtain the short exact sequence 
\begin{equation}\label{kortHfoljd}
0 \to H(l)_{-l} \to H(l)_{-l+1} \to H(l)_{-l+1}/H(l)_{-l} \to 0.
\end{equation}
We are going to show that $H(l)_{-l+1}/H(l)_{-l}\in\add(H)$.
Since $H(l)_{-l}, H(l)_{-l+1} \in \mathcal{F}(\bnabla)$, and
$\mathcal{F}(\bnabla)$ is closed with respect to taking
cokernels of monomorphisms, we obtain that   $H(l)_{-l+1}/H(l)_{-l}
\in \mathcal{F}(\bnabla)$ as well. Hence we can apply the Ringel  
duality functor $F$ to \eqref{kortHfoljd} and obtain the exact sequence 
\begin{equation}\label{kortHfoljd2}
0 \to F(H(l)_{-l}) \to F(H(l)_{-l+1}) \to F(H(l)_{-l+1}/H(l)_{-l}) 
\to 0
\end{equation}
in $R\mathrm{-mod}$. The modules 
$F(H(l)_{-l})$ and $F(H(l)_{-l+1})$ are tilting $R$-modules and
therefore are contained in 
$\mathcal{F}(\bnabla^{(R)})$. Thus  
$F(H(l)_{-l}/H(l)_{-l+1})$ is contained in 
$\mathcal{F}(\bnabla^{(R)})$ as well by the same arguments as above. 
Moreover, the modules $F(H(l)_{-l})$, $F(H(l)_{-l+1})$, and 
$F(H(l)_{-l}/H(l)_{-l+1})$ are contained in 
$\mathcal{F}(\bDelta^{(R)})$ by the Ringel duality. 

Both $F(H(l)_{-l})$ and $F(H(l)_{-l+1})$ have finite projective dimension 
and therefore the projective dimension of $F(H(l)_{-l+1}/H(l)_{-l})$
is also finite. To show that $F(H(l)_{-l+1}/H(l)_{-l})\in 
\mathcal{F}(\Delta^{(R)})$, we use the the following lemma:

\begin{lemma}\label{identdelta}
Let $\mathtt{A}$ be a properly stratified algebra. Then
\begin{enumerate}[(i)] 
\item\label{ident1} $\mathcal{F}(\Delta^{(\mathtt{A})})= \{M^{(\mathtt{A})} 
\in \mathcal{F}(\bDelta^{(\mathtt{A})})\,| \, \text{p.d.}(M^{(\mathtt{A})})
<\infty \}$, and 
\item\label{ident2} $\mathcal{F}(\nabla^{(\mathtt{A})})= \{M^{(\mathtt{A})} 
\in \mathcal{F}(\bnabla^{(\mathtt{A})}) \,| \, \text{i.d.}(M^{(\mathtt{A})})
<\infty \}$.
\end{enumerate}
\end{lemma}

\begin{proof}
We prove \eqref{ident1}. The statement \eqref{ident2} is proved by
similar arguments.

The inclusion $\mathcal{F}(\Delta^{(\mathtt{A})}) \subset \{M^{(\mathtt{A})} 
\in \mathcal{F}(\bDelta^{(\mathtt{A})})\,| \, \text{p.d.}(M^{(\mathtt{A})})
<\infty \}$ can be found for example in \cite[Proposition~1.3]{PR} or in
\cite[Proposition~1.8]{AHLU}. 

Let us prove the inverse inclusion. Let $M^{(\mathtt{A})}
\in \mathcal{F}(\bDelta^{(\mathtt{A})})$ and $\lambda \in \Lambda$.
Consider a short exact sequence,
\begin{equation}\label{nablaseq} 
0 \to  X^{(\mathtt{A})} \to \nabla^{(\mathtt{A})}(\lambda) \to 
\bnabla^{(\mathtt{A})}(\lambda) \to 0,
\end{equation}
where $X^{(\mathtt{A})}$ has a filtration with subquotients 
$\bnabla^{(\mathtt{A})}(\lambda)$. Applying $\Hom_\mathtt{A}(M^{(\mathtt{A})},{}_-)$ 
to \eqref{nablaseq} we get the long exact sequence  
\begin{multline} 
\dots \to \Ext_\mathtt{A}^i(M^{(\mathtt{A})},X^{(\mathtt{A})}) \to \Ext_\mathtt{A}^i(M^{(\mathtt{A})},\nabla^{(\mathtt{A})}(\lambda)) \to 
\Ext_\mathtt{A}^i(M^{(\mathtt{A})},\bnabla^{(\mathtt{A})}(\lambda)) \to \\
\to \Ext_\mathtt{A}^{i+1}(M^{(\mathtt{A})},X^{(\mathtt{A})}) \to \Ext_\mathtt{A}^{i+1}(M^{(\mathtt{A})},\nabla^{(\mathtt{A})}(\lambda)) \to \dots.
\end{multline}
But $\Ext_\mathtt{A}^i(M^{(\mathtt{A})},\nabla^{(\mathtt{A})}(\lambda)) = 0$  
for $i>0$, which gives us a dimension shift between the spaces
$\Ext_\mathtt{A}^i(M^{(\mathtt{A})},\bnabla^{(\mathtt{A})}(\lambda))$ and
$\Ext_\mathtt{A}^{i+1}(M^{(\mathtt{A})},X^{(\mathtt{A})})$. Since  
$\text{p.d.}(M^{(\mathtt{A})})$ is finite, we derive that
$\Ext_\mathtt{A}^i(M^{(\mathtt{A})},\bnabla^{(\mathtt{A})}(\lambda))= 0$ for all  
$i$ big enough. But the dimension shift and the fact that $X^{(\mathtt{A})}$ has a 
filtration with subquotients $\bnabla^{(\mathtt{A})}(\lambda)$ imply that
$\Ext_\mathtt{A}^i(M^{(\mathtt{A})},\bnabla^{(\mathtt{A})}(\lambda))= 0$ for 
all $i$. Hence $M^{(\mathtt{A})} \in \mathcal{F}(\Delta^{(\mathtt{A})})$ by 
\cite[Theorem~1.6(iii)]{AHLU}. This completes the proof of 
Lemma~\ref{identdelta}. 
\end{proof}
 
From Lemma~\ref{identdelta} it follows that $F(H(l)_{-l+1}/H(l)_{-l})$
has a standard filtration and thus is a tilting module. Applying
$F^{-1}$ we obtain  $H(l)_{-l+1}/H(l)_{-l}\in\add(H)$. From 
Lemma~\ref{extHH} it now follows that the short 
exact sequence \eqref{kortHfoljd} splits. This means that $H(l)_{-l}$ 
is a direct summand of $H(l)_{-l+1}$ and by deleting it we construct the 
complex $\mathcal{H}(l-1)^{\bullet}$. 

This proves the existence of the complex $\mathcal{H}(0)^{\bullet}$,
which happens to be a positive complex in $\Ka^b(\add(H))$
quasi-isomorphic to $\mathcal{P}^{\bullet}$. This completes the
proof. 
\end{proof}

\begin{corollary}\label{CC1}
Let $M\in\mathcal{P}(A)^{<\infty}$. Then the module $M$ belongs to
$\add(H)$ if and only if $\Ext_A^i(X,M)=0$ for all $i>0$ and all 
$X\in\mathcal{P}(A)^{<\infty}$.
\end{corollary}

\begin{proof}
Let $M\in\add(H)$. That  $\Ext_A^i(X,M)=0$ for all $i>0$ and all 
$X\in\mathcal{P}(A)^{<\infty}$ follows easily from $\Ext_A^i(H,H)=0$, $i>0$, 
and Proposition~\ref{coresofM} by induction on the length of the
$\add(H)$-coresolution of $M$. 

Let $M\in\mathcal{P}(A)^{<\infty}$ be such that $\Ext_A^i(X,M)=0$ for 
all $i>0$ and all $X\in\mathcal{P}(A)^{<\infty}$, and let
$H_0\in \add(H)$ be such that $M\hookrightarrow H_0$
(such $H_0$ exists by Proposition~\ref{coresofM}). Then the
cokernel of the latter embedding belongs to $\mathcal{P}(A)^{<\infty}$ and 
hence the above condition on $M$ implies that $M$ is in fact a direct
summand of $H_0$, that is $M\in \add(H)$.
\end{proof}

\begin{theorem}\label{thmfdpdh}
Let $A$ be an SSS-algebra such that the Ringel dual $R$ of $A$ is
properly stratified. Then 
\begin{displaymath}
\text{fin.dim}(A)=\text{p.d.}(H).
\end{displaymath}
\end{theorem}

\begin{proof}
If $0\to X\to Y\to Z\to 0$ is an exact sequence then the long exact sequence
implies that 
\begin{equation}\label{eqlltr}
\text{p.d.}(X)  \leq \max\{\text{p.d.}(Y),\text{p.d.}(Z)\}.
\end{equation}
Let $M\in \mathcal{P}(A)^{<\infty}$. Applying \eqref{eqlltr} inductively 
to the $\add(H)$-coresolution of $M$, constructed in Proposition~\ref{coresofM},
one obtains $\text{p.d.}(M)\leq \text{p.d.}(H)$.  This proves the statement
of the theorem.
\end{proof}

Note that Theorem~\ref{thmfdpdh} implies, in particular, that $\text{fin.dim}(A)$
is finite (we have never used the corresponding result from \cite{AHLU2}).

\subsection{Existence of $H$ guarantees that the category of modules
of finite projective dimension is contravariantly finite}\label{s6.3} 

For a full subcategory, $\mathcal{C}$, of $A$-mod, denote by 
$\check{\mathcal{C}}$ the full subcategory of $A$-mod, which contains
all  modules $M$ for which there is a finite exact sequence, 
\begin{displaymath}
0 \to M \to C_0 \to C_1 \to \dots
\to C_k \to 0,
\end{displaymath}
with $C_i\in\mathcal{C}$. The category
$\hat{\mathcal{C}}$ is defined dually. 

Recall that a full subcategory, $\mathcal{C}$, of $A$-mod is called 
\textit{contravariantly finite} provided that it is closed under
direct summands and isomorphisms, and for each $A$-module $X$ there 
exists a homomorphism $f:C_X \to X$, where $C_X\in\mathcal{C}$, such 
that for any homomorphism $g:C\to X$ with $C\in\mathcal{C}$ there is 
a homomorphism $h:C\to C_X$ such that $f \circ h = g$.

Recall also that a subcategory, $\mathcal{B}$, of $A$-mod is called 
\textit{resolving} if it contains all projective modules and is 
closed under extensions and kernels of epimorphisms.
Obviously, $\mathcal{P}(A)^{<\infty}$ is a resolving category.  
However, $\mathcal{P}^{<\infty}$ is not contravariantly finite in 
general, see \cite{IST}. Our main result in this subsection is the
following:

\begin{theorem}\label{t6.3.1}
Let $A$ be an SSS-algebra, whose Ringel dual is properly
stratified. Then $\mathcal{P}^{<\infty}$ is contravariantly finite
in $A\mathrm{-mod}$. 
\end{theorem}

\begin{proof}
Since $H$ is a (generalized) tilting module, the subcategory 
$\check{\add(H)}$ is contravariantly finite and resolving by 
\cite[Section~5]{AR}. From Proposition~\ref{coresofM} we see that 
$\mathcal{P}^{<\infty}\subset\check{\add(H)}$. On the other hand, 
$\text{p.d.}(H)<\infty$ implies $\mathcal{P}^{<\infty}\supset
\check{\add(H)}$. This completes the proof.
\end{proof}

\subsection{$H$ and $\text{fin.dim.}(A)$ for algebras 
with duality}\label{s6.4} 

In this section we calculate the finitistic dimension of a properly 
stratified algebra $A$ having a simple preserving duality in terms 
of the projective dimension of the characteristic tilting module. 
This generalizes the main result in \cite{MO}, where analogous result
is obtained under the assumption that the characteristic tilting 
and cotilting $A$-modules coincide. The main result of the section is: 

\begin{theorem}\label{findim}
Let $A$ be a properly stratified algebra having a simple preserving 
duality, whose Ringel dual $R$ is also properly stratified. Then
\begin{displaymath}
\text{fin.dim}(A)=2\text{p.d.}(T^{(R)}).  
\end{displaymath}
\end{theorem}

To prove the statement we will need several lemmas.

\begin{lemma}\label{uppskatta}
Let $A$ be as in Theorem~\ref{findim}. Then
\begin{displaymath}
\text{p.d.}(T^{(R)})\leq\text{p.d.}(T).  
\end{displaymath}
\end{lemma}

\begin{proof}
Since $T$ is a tilting module, we can choose a minimal tilting
coresolution,
\begin{displaymath}
0\to P\to T_0\to T_1\to \dots\to T_a\to 0,
\end{displaymath}
of $P$, where $a=\text{p.d.}(T)$. From this coresolution 
we obtain the complex 
\begin{displaymath}
\dots \to 0 \to T_0 \to T_1 \to \dots \to T_a \to 0 \to \dots,
\end{displaymath}
in which the only non-zero homology is in degree zero and equals
$P$. 

On the other hand, the module $F(H)$ is tilting over $R$, in particular,
it has finite projective dimension. Hence 
we can choose a minimal projective resolution,
\begin{displaymath}
0 \to  P^{(R)}_b \to \dots  \to P^{(R)}_1 \to  
P^{(R)}_0 \to F(H) \to 0,
\end{displaymath}
of $F(H)$, where $b=\text{p.d.}(T^{(R)})$.  Applying $F^{-1}$ 
we get a minimal tilting resolution,
\begin{equation}\label{eqref2}
0 \to T_{(b)}=F^{-1}(P^{(R)}_b) \to \dots \to T_{(1)}=F^{-1}(P^{(R)}_1) \to 
T_{(0)}=F^{-1}(P^{(R)}_0) 
\to H  \to 0,
\end{equation}
of $H$. Hence we also obtain the complex   
\begin{displaymath}
\mathcal{J}^\bullet:\quad\quad \dots \to 0 \to T_{(b)} \to
\dots \to T_{(1)} \to T_{(0)} \to 0 \to \dots,
\end{displaymath}
in which the only non-zero homology is in degree zero and equals 
$H$.

If $0\to X\to Y\to Z\to 0$ is an exact sequence, then
$\mathrm{p.d.}(Z)\leq \max\{\mathrm{p.d.}(X),\mathrm{p.d.}(Y)\}+1$. 
Applying this inequality inductively to \eqref{eqref2} we get
\begin{equation}\label{eq9.22}
\text{p.d.}(H)\leq a+b.
\end{equation}

Applying ${}^\circ$ to $\mathcal{J}^\bullet$ gives the complex 
\begin{displaymath}
\mathcal{C}^\bullet:\quad\quad \dots \to 0 \to 
T_{(0)}^\circ \to
T_{(1)}^\circ \to \dots \to 
T_{(b)}^\circ\to 0 \to \dots,
\end{displaymath}
where the only non-zero homology is in degree zero and equals 
$H^\circ$. Remark that the complex
$\mathcal{C}^\bullet$ consists of cotilting $A$-modules.
We would like to substitute $\mathcal{C}^\bullet[b]$ 
by a negative complex from $\Ka(\add(T))$. To be able to do this 
we need the following lemma:

\begin{lemma}\label{tiltresofC}
Let $A$ be a properly stratified algebra. Then for each
$\lambda\in\Lambda$ there exists a (possibly infinite)
minimal tilting resolution of $C(\lambda)$, which has the following 
form:
\begin{displaymath}
 \dots \to T_1 \to T(\lambda)\oplus T_0 \to C(\lambda) \to 0.
\end{displaymath}
\end{lemma}

\begin{proof}
Let $\lambda\in\Lambda$. Since $C(\lambda)$ has a proper
costandard filtration, we can apply the Ringel duality functor $F$ and 
get the $R$-module $F(C(\lambda))$. Choose a minimal projective 
resolution, 
\begin{equation} \label{hoplya}
\dots \to P^{(R)}_1 \to P^{(R)}_0 \to F(C(\lambda)) \to 0, 
\end{equation}
of $F(C(\lambda))$. Since $C(\lambda)$ surjects onto 
$\nabla(\lambda)$ and $\nabla(\lambda)$ has a filtration with 
subquotients $\bnabla(\lambda)$, 
applying $F$ it follows that there is an epimorphism from 
$F(C(\lambda))$ to $\bDelta^{(R)}(\lambda)$. 
Hence we conclude that the head of 
$F(C(\lambda))$ contains $L^{(R)}(\lambda)$, and therefore
$P^{(R)}_0=P^{(R)}(\lambda)\oplus \hat P^{(R)}$, where $\hat P^{(R)}$ is 
some projective $R$-module. Since $F(C(\lambda))\in \mathcal{F}(\bDelta^{(R)})$
and \eqref{hoplya} is exact, it follows that the kernels of all morphisms in 
\eqref{hoplya} belong to $\mathcal{F}(\bDelta^{(R)})$ as well. The statement 
of the lemma now follows by applying $F^{-1}$ to \eqref{hoplya} and
taking $P^{(R)}_0=P^{(R)}(\lambda)\oplus \hat P^{(R)}$ into account.
\end{proof}

From Lemma~\ref{tiltresofC} it follows that we can apply
Lemma~\ref{apl2} to $\mathtt{A}=A$, $X^{(\mathtt{A})}=C$, 
$Y^{(\mathtt{A})}=T$ and the complex 
$\mathcal{C}^\bullet[b]$, and obtain a negative complex in
$\Ka(\add(T))$, which is quasi-isomorphic to $\mathcal{C}^\bullet[b]$.
Shifting the latter one by $-b$ in the derived category gives
a complex, 
\begin{displaymath}
\mathcal{T}^\bullet:\quad\quad \dots \to 
\ T(0)\to\dots\to T(b-1)\to T(b)\to 0\to \dots,
\end{displaymath}
in $\Ka^-(\add(T))$ with the only non-zero homology in degree
zero, which equals $H^\circ$. We assume that $\mathcal{T}^\bullet$
is minimal that is does not contain trivial direct summands.
Lemma~\ref{tiltresofC} implies that $T(b)=T_{(b)}\oplus
\hat{T}$ for some tilting module $\hat{T}$.

We have
\begin{displaymath}
\Ext_A^{2b}\left(H,H^\circ\right)=
\Hom_{\Der^-(A)}\left(H^{\bullet}[-b],\left(H^\circ\right)^{\bullet}[b]\right)=
\Hom_{\Der^-(A)}\left(\mathcal{J}^\bullet[-b],
\mathcal{T}^\bullet[b]\right).
\end{displaymath}
From \cite[Chap.~III, Lemma~2.1]{Ha} we have 
\begin{displaymath}
\Hom_{\Der^-(A)}\left(\mathcal{J}^\bullet[-b],
\mathcal{T}^\bullet[b]\right)=
\Hom_{\Ka^-(A)}\left(\mathcal{J}^\bullet[-b],
\mathcal{T}^\bullet[b]\right).  
\end{displaymath}
Let $f:T_{(b)} \hookrightarrow T(b)$ be the inclusion, defined via the isomorphism
of $T_{(b)}$  with the first direct summand of $T(b)$. Denote by $g:T_{(b)}\to 
T_{(b-1)}$ and $h:T(b-1)\to T(b)$ the differentials
in the complexes $\mathcal{J}^\bullet$ and
$\mathcal{T}^\bullet$ respectively and consider the following diagram:
\begin{displaymath}
\xymatrix{
 & F\left(T_{(b)}\right)\ar[d]^{F(f)}\ar[r]^{F(g)}\ar@{.>}[dl]_{\alpha} 
&F\left(T_{(b-1)}\right)\ar@{.>}[dl]_{\beta} \\
F\left(T(b-1)\right)\ar[r]^{F(h)} & F\left(T(b)\right)& 
}
\end{displaymath}
The minimality of
$\mathcal{J}^\bullet[-b]$ and $\mathcal{T}^\bullet[b]$ implies
that the images of the morphisms $F(g)$ and $F(h)$ belong to the
radicals of the corresponding modules. Hence for every $\alpha$ and
$\beta$ as depicted on the diagram the image of
$F(h)\circ\alpha+\beta\circ F(g)$ belongs to the radical of
$F(T(b))$. However the image of $F(f)$ does not belong to the
radical of $F(T(b))$. This means that the morphism $f$ induces a
non-zero homomorphism from $\mathcal{J}^\bullet[-b]$ to 
$\mathcal{T}^\bullet[b]$ in $K^-(A)$.   
From this we conclude that
$\Ext_A^{2b}\left(H,H^\circ\right)\neq 0$
and hence 
\begin{equation}\label{eq9.23}
2b \leq \text{p.d.}(H).
\end{equation}
Combining \eqref{eq9.22} and \eqref{eq9.23} we obtain $b\leq a$, which
completes the proof of Lemma~\ref{uppskatta}. 
\end{proof}

The arguments above immediately imply:

\begin{corollary}\label{c6.3.22}
$\Ext_A^{i}\left(H,H^\circ\right)=0$ for all 
$i>2\text{p.d.}(T^{(R)})$.
\end{corollary}

Further, we can derive the following inequality
(compare with \cite[Theorem~2.2.1]{EP}):

\begin{lemma}\label{uppskatta2}
Let $A$ be as in Theorem~\ref{findim} and $k = \text{p.d.}(H)$. Then
$\Ext^k_A\left(H,H^\circ\right)\neq 0$.
\end{lemma}

\begin{proof}
Choose a minimal projective resolution, 
\begin{displaymath}
0\to P_k\to\dots\to P_1 \to P_0 \to H\to 0, 
\end{displaymath}
of $H$ and let 
\begin{displaymath}
\mathcal{P}^\bullet:\quad\quad
\dots\to 0\to P_k\to\dots\to P_1\to P_0 \to 0 \to \dots 
\end{displaymath}
be the corresponding complex in $\Ka^b(\add(P))$. 
Choose also a minimal (possibly infinite) 
projective resolution,   
\begin{displaymath}
\dots\to Q_1 \to Q_0 \to H^\circ \to 0,
\end{displaymath}
of $H^\circ$, and construct the corresponding (possibly
infinite) complex  
\begin{displaymath}
\mathcal{Q}^\bullet:\quad\quad \dots \to Q_1 \to Q_0 \to 0 \to
\dots. 
\end{displaymath}
Applying ${}^\circ$ to the short exact sequence 
\begin{displaymath}
0 \to \bnabla \to H \to \text{Coker} \to 0
\end{displaymath}
given by Proposition~\ref{prref3}, we obtain the short exact sequence  
\begin{displaymath}
0 \to \text{Coker}^\circ \to H^\circ \to 
\bDelta \to 0.
\end{displaymath}
It follows that
the head of $H^\circ$ contains the head of
$\bDelta$, which coincides with $L$. 
Hence $Q_0=P\oplus Q$, where $Q$ is
some projective module. Using the same arguments as in the proof of
Lemma~\ref{uppskatta} we obtain 
$\Hom_{\Der^-(A)}(\mathcal{P}^\bullet,\mathcal{Q}^\bullet[k])\neq 0$
and hence $\Ext_A^{k}\left(H,H^\circ\right)\neq 0$.
This completes the proof.
\end{proof}

Now we are ready to prove Theorem~\ref{findim}.

\begin{proof}[Proof of Theorem~\ref{findim}.]
Let $k=\text{p.d.}(H)$. 
Using Theorem~\ref{thmfdpdh} we have
$k=\text{fin.dim}(A)$.
From Lemma~\ref{uppskatta2} we obtain that 
$\Ext^k_A\left(H,H^\circ\right)\neq
0$. Certainly, 
$\Ext^i_A\left(H,H^\circ\right)= 0$ for all
$i>k$. But from the proof of Lemma~\ref{uppskatta} and from
Corollary~\ref{c6.3.22} we conclude $k=2b$. This completes the proof.
\end{proof}

To relate the finitistic dimension of $A$ to $\mathrm{p.d.}(T)$ 
we will need a stronger assumption.

\begin{proposition}\label{c63.25}
Let $A$ be as in Theorem~\ref{findim} and assume that $R$ also
has a  simple preserving duality.  Then 
\begin{displaymath}
\text{fin.dim}(A)=2\text{p.d.}(T).  
\end{displaymath}
\end{proposition}

\begin{proof}
From Lemma~\ref{uppskatta} we have
$\text{p.d.}(T^{(R)})\leq\text{p.d.}(T)$. The existence of the
duality for both $A$ and $R$ implies that 
$\text{p.d.}(T^{(R)})=\text{i.d.}(C^{(R)})$ and
$\text{p.d.}(T)=\text{i.d.}(C)$.
Applying Lemma~\ref{uppskatta} to $A^\mathrm{opp}$ and $R^\mathrm{opp}$
and using the usual duality we obtain 
$\text{i.d.}(C)\leq\text{i.d.}(C^{(R)})$.
Altogether we deduce $\text{p.d.}(T^{(R)})=\text{p.d.}(T)$ and
the statement follows from Theorem~\ref{findim}.
\end{proof}

We remark that if $A$ is as in Theorem~\ref{findim}, then $R$ does not
necessarily have a simple preserving duality, see the example 
in Subsection~\ref{s8.25}.   

\section{Two-step duality for standardly stratified
algebras}\label{s55}  

Since $H$ is a (generalized) tilting module, the classical tilting theory
suggests to study the algebra $B(A) = \text{End}_{A}(H)$. Consider the 
functor $G:A\text{-mod}\to B(A)\text{-mod}$ defined via 
\begin{displaymath}
G({}_-) = \DD \circ \Hom_A({}_-,H).
\end{displaymath}  
From the definition it follows that $G(H)$ is an injective cogenerator of
$B(A)$. 

Consider also the functor $G':B(A)\text{-mod}\to A\text{-mod}$ defined 
via 
\begin{displaymath}
G'({}_-) = \Hom_{B(A)}(G(A),{}_-).
\end{displaymath}  

We start with establishing some necessary properties of the functor $G$.

\begin{lemma}\label{l01}
\begin{enumerate}
\item $G$ is exact on $\mathcal{P}(A)^{<\infty}$, and maps 
$\mathcal{P}(A)^{<\infty}$ to $\mathcal{I}(B(A))^{<\infty}$.
\item $G$ is full and faithful on $\mathcal{P}(A)^{<\infty}$.
\end{enumerate}
\end{lemma}

\begin{proof}
That $G$ is exact on $\mathcal{P}(A)^{<\infty}$ follows from
Corollary~\ref{CC1}. 
Let $M\in A\text{-mod}$ be such that $\text{p.d.}(M)<\infty$. By 
Proposition~\ref{coresofM} there exists a coresolution 
\begin{equation}\label{FF1F}
0 \to M \to H_0 \to \dots \to H_k \to 0,
\end{equation}
where $H_i \in \add(H)$ and $k\geq 0$. Applying the exact functor $G$
gives the exact sequence 
\begin{displaymath}
0 \to G(M) \to G(H_0) \to \dots \to G(H_k) \to 0,
\end{displaymath}
in $B(A)$-mod. Since $G(H_j)$ is $B(A)$-injective for all $j$, we obtain 
that the injective dimension of $G(M)$ is finite. This proves the first
statement. 

Let us now prove that $G$ is full and faithful on
$\mathcal{P}(A)^{<\infty}$. We start with showing that $G$ 
is full and faithful on $\add(H)$.
For this we calculate the following: 
\begin{multline}\label{FFF1}
G'\circ G({}_A H_{B(A)}) = \\ = \Hom_{B(A)-}\left(\DD\circ 
\Hom_A({}_A A_{A},{}_A H_{B(A)}),\DD\circ 
\Hom_A({}_A H_{B(A)},{}_A H_{B(A)})\right)= \\ =  
\Hom_{-B(A)}\left(\Hom_A({}_A H_{B(A)},{}_A H_{B(A)}),
\Hom_A({}_A A_{A},{}_A H_{B(A)})\right) = \\ =
\Hom_{-B(A)}\left(B(A),{}_A H_{B(A)}\right)= {}_A H_{B(A)},
\end{multline}
and 
\begin{multline}\label{FFF2}
G\circ G'({}_{B(A)} \DD(B(A)^\mathrm{opp})) = \\ =
\DD\circ \Hom_A\left(\Hom_{B(A)-}
\left(\DD\circ \Hom_A({}_A A_{A},{}_A H_{B(A)}),{}_{B(A)} 
\DD(B(A)^\mathrm{opp})\right),{}_A H_{B(A)}\right)= \\ =
\DD\circ \Hom_A\left(\Hom_{-B(A)}
\left(B(A)_{B(A)},\Hom_A({}_A A_{A},{}_A H_{B(A)})\right),
{}_A H_{B(A)}\right) =\\ =
\DD\circ \Hom_A\left(\Hom_{-B(A)}
\left(B(A)_{B(A)},{}_A H_{B(A)}\right),
{}_A H_{B(A)}\right) = \\ =
\DD\circ \Hom_A\left({}_A H_{B(A)},{}_A H_{B(A)}\right) 
= {}_{B(A)} \DD(B(A)^\mathrm{opp}),
\end{multline}
which implies that $G$ is full and faithful on $\add(H)$.

Now the fact that $G$ is full and faithful on $\mathcal{P}(A)^{<\infty}$ 
follows from the existence of \eqref{FF1F} for all 
$M\in \mathcal{P}(A)^{<\infty}$ by induction on the length of the
coresolution \eqref{FF1F}. This completes the proof.
\end{proof}

\begin{lemma}\label{l015}
\begin{enumerate}
\item The functor $G$ maps $N(\lambda)$ to
$\nabla^{(B(A))}(\lambda)$ for all $\lambda\in\Lambda$.
\item The algebra $(B(A)^{\mathrm{opp}},\leq)$ is an SSS-algebra.  
\end{enumerate}
\end{lemma}

\begin{proof}
For $\lambda \in \Lambda$ consider the short exact sequence 
\begin{displaymath}
0\to N(\lambda)\to H(\lambda)\to Y(\lambda)\to 0,
\end{displaymath} 
where $Y(\lambda)$ has a filtration with subquotients
$N(\mu)$, $\mu>\lambda$ (see Proposition~\ref{prref3}). 
Because of the exactness of $G$ on 
$\mathcal{P}(A)^{<\infty}$, see 
Lemma~\ref{l01}, and the fact that $\text{p.d.}(N)<\infty$,
the sequence above yields the short exact sequence 
\begin{displaymath}
0 \to G\left(N(\lambda)\right) \to I^{(B(A))}(\lambda) 
\to G\left(Y(\lambda)\right) \to 0,
\end{displaymath} 
where $G\left(Y(\lambda)\right)$ has a filtration with subquotients
$G\left(N(\mu\right))$, $\mu>\lambda$. Further, using the
classical Ringel duality, for $\lambda<\mu$ we have 
\begin{displaymath}
\Hom_A\left(N(\lambda),H(\mu)\right) =
\Hom_R\left(\Delta^{(R)}(\lambda),T^{(R)}(\mu)\right) = 0
\end{displaymath}
(we recall that $\Delta^{(R)}(\lambda)$ are defined with respect to 
$\leq_{R}$, which is opposite to the original order $\leq$).
Using the fact that $G$ is full and faithful on 
$\mathcal{F}(N)\subset \mathcal{P}(A)^{<\infty}$ we obtain  
\begin{multline*}
\left[G\left(N(\lambda)\right):L^{(B(A))}(\mu)\right] = 
\dim_\k\Hom_{B(A)}\left(G\left(N(\lambda)\right),
I^{(B(A))}(\mu)\right) = \\
= \dim_\k\Hom_{A}\left(N(\lambda),H(\mu)\right) = 0.  
\end{multline*}

So the family $\{G\left(N(\lambda)\right)|\lambda \in \Lambda \} 
\subset B(A)\text{-mod}$ satisfies all the conditions, 
characterizing the costandard modules for SSS-algebras
(see for example \cite[Lemma~1.5]{AHLU}). This 
implies $G\left(N(\lambda)\right) = \nabla^{(B(A))}(\lambda)$. 
Exactness of 
$G$ on $\mathcal{F}(N)\subset \mathcal{P}(A)^{<\infty}$ 
guarantees that the injective cogenerator $G(H)$  of $B(A)$ is filtered by 
costandard modules. This completes the proof of both statements of 
the lemma.
\end{proof}

\begin{lemma}\label{l02}
\begin{enumerate}
\item $G(A)$ is a (generalized) cotilting module for $B(A)$.
\item $G(A)\in\mathcal{F}\left(\bDelta^{(B(A))}\right)$. 
\end{enumerate}
\end{lemma}

\begin{proof}
Using Lemma~\ref{l01} one easily obtains 
\begin{displaymath}
\Ext_{B(A)}^i\big(G(A),G(A)\big)=\Ext_{A}^i\big(A,A\big) = 0
\end{displaymath}
for all $i>0$. Since $\text{p.d.}(A)=0<\infty$, from  Lemma~\ref{l01} 
it also follows that $\text{i.d.}(G(A))<\infty$. Since $H$ has
finite projective dimension, we can take a minimal projective 
resolution, 
\begin{displaymath}
0 \to P_k \to \dots  \to P_1 \to P_0 \to H \to 0,
\end{displaymath}
of $H$.
Since all modules in this resolution have finite projective dimension 
we can apply $G$  and use Lemma~\ref{l01} to obtain a resolution of 
the injective cogenerator $G(H)$ over $B(A)$ by modules $G(P_i)$ from  
$\add(G(A))$. Hence $G(A)$ is a (generalized) cotilting module for
$B(A)$. 

Using Lemma~\ref{l01} and Lemma~\ref{l015} one gets 
\begin{displaymath}
\Ext_{B(A)}^i\left(G(A),\nabla^{(B(A))}\right)=
\Ext_{B(A)}^i\left(G(A),G\left(N\right)\right)=
\Ext_{A}^i\left(A,N\right) = 0,
\end{displaymath}
implying $G(A)\in\mathcal{F}\left(\bDelta^{(B(A))}\right)$. This completes the 
proof.
\end{proof}

\begin{lemma}\label{l03}
The functor $G$ induces an equivalence between $\mathcal{F}(N)$ 
and $\mathcal{F}(\nabla^{(B(A))})$ with the inverse functor $G'$.
\end{lemma}

\begin{proof}
By  Lemma~\ref{l01} for $M\in\mathcal{P}(A)^{<\infty}$ we have 
\begin{equation}\label{FF2F}
G'\circ G(M)=\Hom_{B(A)}(G(A),G(M))= \Hom_{A}(A,M)= M,
\end{equation}
which implies that $G'\circ G$ is isomorphic to the identity functor 
on $\mathcal{F}(N)$.

After Lemma~\ref{l01} we have that $G$ is full and faithful on 
$\mathcal{F}(N)$, and thus we have only to prove that it is
dense.  

From the second statement of Lemma~\ref{l02} we obtain that $G'$ is
exact on $\mathcal{F}(\nabla^{(B(A))})$. \eqref{FFF1} and \eqref{FFF2} 
imply that 
$G'$ is full and faithful on $\add(\DD(B(A)^\mathrm{opp}))$ that is on all
$B(A)$-injective modules. Since every module with a costandard
filtration has a finite coresolution by injective modules we obtain
that $G'$ is full and faithful on $\mathcal{F}(\nabla^{(B(A))})$. But 
\eqref{FF2F} implies that 
$G':\mathcal{F}(\nabla^{(B(A))})\to\mathcal{F}(N)$ and 
 is dense. Thus $G'$ is an equivalence and \eqref{FF2F} implies that $G$
is inverse to $G'$, hence is dense as well. This completes the proof.
\end{proof}

Now we can formulate the main result of this section.

\begin{theorem}\label{thmtsd}
Let $(A,\leq)$ be an SSS-algebra, whose Ringel dual $(R,\leq_{R})$ is 
properly stratified. Then 
\begin{enumerate}[(i)] 
\item\label{2Step1} The algebra $B(A)^{\mathrm{opp}}$ is an SSS-algebra
and is isomorphic to the opposite algebra of the second Ringel 
dual $End_{R}(T^{(R)})$. 
\item\label{2Step2} $B(A)^\mathrm{opp}$ has the Ringel dual
$(R^\mathrm{opp},\leq_{R})$, which is properly stratified, 
and the algebra $B(B(A)^\mathrm{opp})^{\mathrm{opp}}$ is Morita 
equivalent to $A$.   
\end{enumerate}
\end{theorem}

\begin{proof} 
That $(B(A)^{\mathrm{opp}},\leq)$ is an SSS-algebra was proved 
in Lemma~\ref{l015}. The statement about the second Ringel dual
follows from the usual Ringel duality (\cite[Theorem~2.6]{AHLU}) by
\begin{displaymath}
\End_A(H) = \End_R\left(F\big(H\big)\right) = \End_R(T^{(R)}),
\end{displaymath}
since $F\big(H\big) = T^{(R)}$. This proves \eqref{2Step1}. 

Now we prove \eqref{2Step2}. We start with calculation of the 
indecomposable tilting modules in $B(A)^{\mathrm{opp}}$. Composing 
the functor $G$ with the duality $\DD:B(A)\text{-mod} \to 
B(A)^{\mathrm{opp}}\text{-mod}$, we get the contravariant functor 
$\DD \circ G:A\text{-mod}\to B(A)^{\mathrm{opp}}\text{-mod}$. Applying 
this functor to the short exact sequence \eqref{seq1} and using the 
exactness of  $G$ on  $\mathcal{F}(N)$, we get the exact sequence 
\begin{displaymath}
0 \to \DD\circ G(N(\lambda)) \to \DD\circ G(T(\lambda))\to
\DD\circ G(S(\lambda)) \to 0,
\end{displaymath}
where $\DD\circ G(N(\lambda)) = \DD \nabla^{(B(A))}(\lambda) = 
\Delta^{(B(A)^{\mathrm{opp}})}(\lambda)$, and $\DD\circ 
G(S(\lambda))$ has a filtration with subquotients   
$\Delta^{(B(A)^{\mathrm{opp}})}(\mu)$, $\mu < \lambda$. Moreover, 
using Lemma~\ref{l03}, we have 
\begin{displaymath}
0 = \Ext_A^i\left(T(\lambda),N(\mu)\right) \cong 
\Ext_{B(A)^{\mathrm{opp}}}^i\left(\Delta^{(B(A)^{\mathrm{opp}})}(\mu),
\DD\circ G\big(T(\lambda)\big)\right)
\end{displaymath}
for all $\lambda$, $\mu$ and $i>0$. Hence $\DD\circ G\big(T(\lambda)\big) \in 
\mathcal{F}(\bnabla^{(B(A)^{\mathrm{opp}})})$
and we can conclude that $\DD\circ G\big(T(\lambda)\big)$ is an
indecomposable tilting module, and, moreover, that
$\DD\circ G\big(T(\lambda)\big) = T^{(B(A)^{\mathrm{opp}})}(\lambda)$.
The Ringel dual to $B(A)^{\mathrm{opp}}$ is now computed by 
\begin{displaymath}
\End_{B(A)^{\mathrm{opp}}}\left(T^{(B(A)^{\mathrm{opp}})}\right) =
\End_A\left(T\right)^{\mathrm{opp}} = R^{\mathrm{opp}},  
\end{displaymath}
and thus is properly stratified.

Hence we have the corresponding functor
$F^{B(A)^{\mathrm{opp}}}:B(A)^{\mathrm{opp}}\text{-mod} \to 
R^{\mathrm{opp}}\text{-mod}$. Since $R^{\mathrm{opp}}$ is also 
properly stratified, we can construct $N^{(B(A)^{\mathrm{opp}})}(\lambda)=  
(F^{B(A)^{\mathrm{opp}}})^{-1}\left(\Delta^{(R^{\mathrm{opp}})}(\lambda)\right)$
and
$H^{(B(A)^{\mathrm{opp}})}(\lambda) = 
(F^{B(A)^{\mathrm{opp}}})^{-1}T^{(R^{\mathrm{opp}})}(\lambda)$. Then
$B(B(A)^{\mathrm{opp}}) = 
\End_{B(A)^{\mathrm{opp}}}(H^{B(A)^{\mathrm{opp}}})$ and 
we compute  
\begin{multline*}
B(B(A)^{\mathrm{opp}})^{\mathrm{opp}} = 
\left(\End_{B(A)^{\mathrm{opp}}}
\big(H^{(B(A)^{\mathrm{opp}})}\big)\right)^{\mathrm{opp}}  = 
\left(\End_{R^{\mathrm{opp}}}
\big(T^{(R^{\mathrm{opp}})}\big)\right)^{\mathrm{opp}} = \\
= \End_{R}\big(C^{(R)}\big) = 
\End_{A}\big(I\big) \eqsim A, 
\end{multline*}
where $\eqsim$ denotes Morita equivalence. This completes the proof. 
\end{proof}

Lemma~\ref{l01} and Lemma~\ref{l03} admit the following extension:

\begin{proposition}\label{p04}
The functor $G$ induces an equivalence between $\mathcal{P}(A)^{<\infty}$
and  $\mathcal{I}(B(A))^{<\infty}$ with the inverse functor $G'$.
\end{proposition}

To prove this we will need the following lemma:

\begin{lemma}\label{l05}
Let $M^{(B(A))}\in\mathcal{I}(B(A))^{<\infty}$. Then there exists a resolution, 
\begin{displaymath}
0 \to Y_k^{(B(A))} \to \dots \to Y_1^{(B(A))} \to Y_0^{(B(A))} \to M^{(B(A))} \to 0,
\end{displaymath}
where $Y_i^{(B(A))}\in\add(G(A))$ for all $i$. 
\end{lemma}

\begin{proof}
Consider an injective coresolution,  
\begin{displaymath}
0 \to M^{(B(A))} \to I_0^{(B(A))} \to I_1^{(B(A))} \to \dots \to I_m^{(B(A))} \to 0,
\end{displaymath}
of $M^{(B(A))}$ and let $\mathcal{M}^{\bullet}$ be the corresponding complex in  
$\Ka(I^{(B(A))})$. Applying Lemma~\ref{apl1} to
$\mathtt{A}=B(A)^{\mathrm{opp}}$, $X^{(\mathtt{A})}=\DD(I^{(B(A))})$, 
$Y^{(\mathtt{A})}=\DD(G(A))$ and
the complex $\DD(\mathcal{M}^{\bullet})$ gives a finite complex, the dual
of which belongs to $\Ka(\add(G(A)))$ and is quasi-isomorphic to
$\mathcal{M}^{\bullet}$. The necessary resolution is now obtained using 
the projectivity of ${}_A A$, properties of $G$ given by
Lemma~\ref{l01}, and arguments dual to those used in
Proposition~\ref{coresofM}. We omit the details. 
\end{proof}

We are now ready to prove Proposition~\ref{p04}.

\begin{proof}[Proof of Proposition~\ref{p04}.]
By  Lemma~\ref{l01} for $M\in\mathcal{P}(A)^{<\infty}$ we have 
\begin{equation}\label{FF2}
G'\circ G(M)=\Hom_{B(A)}(G(A),G(M))= \Hom_{A}(A,M)= M,
\end{equation}
which implies that $G'\circ G$ is isomorphic to the identity functor
on $\mathcal{P}(A)^{<\infty}$. From Lemma~\ref{l01} we have that $G$ 
is full and faithful on $\mathcal{P}(A)^{<\infty}$ and thus we
are only left to prove that it is dense. From Lemma~\ref{l05}, using 
induction on the length of the $\add(G(A))$-resolution, we obtain that 
$\Ext_{B(A)}^i(G(A),X)=0$ for all $i>0$ and $X\in\mathcal{I}(B(A))^{<\infty}$, 
Thus $G'$ is exact on $\mathcal{I}(B(A))^{<\infty}$. One now completes the proof
using the same arguments as in Lemma~\ref{l03}.
\end{proof}

After Theorem~\ref{thmtsd} we can define $N^*=\DD(N^{(B(A)^{\mathrm{opp}})})$
and $H^*=\DD(H^{(B(A)^{\mathrm{opp}})})$. With this notation we have the
following images of the two-step duality functor $G$.

\begin{proposition}\label{p07} 
For every $\lambda \in \Lambda$ we have 
\begin{gather*}
G(H(\lambda))= I^{(B(A))}(\lambda), \quad
G(N(\lambda))= \nabla^{(B(A))}(\lambda), \quad
G(T(\lambda))= C^{(B(A))}(\lambda)\\
G(\Delta(\lambda))= (N^*)^{(B(A))}(\lambda),\quad\quad
G(P(\lambda))= (H^*)^{(B(A))}(\lambda).
\end{gather*}
\end{proposition} 

\begin{proof}
The first three equalities were proved during the proof of
Theorem~\ref{thmtsd}. The fourth equality follows from the third one and the 
fact that $G$ is exact on $\mathcal{P}(A)^{<\infty}$ by induction on $\lambda$,
which starts  from that $\lambda$ for which $T(\lambda)=\Delta(\lambda)$.
The fifth equality follows from the fourth one and the fact that $G$
is exact on $\mathcal{P}(A)^{<\infty}$ by induction on $\lambda$,
which starts  from that $\lambda$ for which $P(\lambda)=\Delta(\lambda)$.
\end{proof}

We would like to end this section with the following two remarks:
firstly, the Ringel and the two-step dualities give rise to the following 
schematic picture of functors on module categories of SSS-algebras:
\begin{displaymath}
\xymatrix@!=0.1pc{
  &    &  R\mathrm{-mod}\ar@{->}[ddrr]^{\DD\circ F^{R}}   & & \\
  &    &      &              & \\
A\mathrm{-mod}\ar@{->}[uurr]^{F}
\ar@{->}[rrrr]^>>>>>>>{\DD\circ G}  & & &    
& B(A)^\text{opp}\mathrm{-mod}\ar@{->}[ddll]^{F^{B(A)^\text{opp}}}\\
  &    &      &              & \\
  &    &  R^\text{opp}\mathrm{-mod}\ar@{->}
[uull]^{\DD\circ F^{R^\text{opp}}}   &    & \\
}
\end{displaymath}
Note that the picture above is {\em not} a commutative diagram. In particular, 
the two-step duality functor $G$ is {\em not} the composition of the Ringel 
dualities $F$ and $F^{R}$. Each functor on this picture induces an equivalence 
of certain subcategories. However, these subcategories are not well-coordinated
with each other. A deeper understanding of the picture above might be
an interesting problem.

Secondly, in Section~\ref{s3} we have shown that the information about the
proper stratification of the Ringel dual of an SSS-algebra, $\mathtt{A}$, can
be derived directly from $\mathtt{A}\mathrm{-mod}$. An interesting problem 
seems to be whether $\mathtt{A}\mathrm{-mod}$ contains enough information
to answer directly the question about the proper stratification of the the two-step dual.

\section{$\mathcal{F}(N)$-filtration dimension}\label{s7}

Recall (see for example \cite[Section~4.2]{MO}) that for an algebra, $A$, a 
family, $\mathcal{M}$, of $A$-modules, and an $A$-module, $M$, one says that the
\textit{$\mathcal{M}$-filtration codimension}
$\codim_{\mathcal{M}}(M)$ of $M$ equals $n$ provided that there 
exists an exact sequence,
\begin{displaymath}
0 \to M \to M_0 \to M_1 \to \dots \to M_n \to 0,
\end{displaymath}
where $M_i \in \mathcal{M}$, and $n$ is minimal with this property.
In this section we study the $\mathcal{F}(N)$-filtration 
codimension for $A$-modules. We start with determining the modules for 
which the notion of $\mathcal{F}(N)$-filtration codimension makes sense.

\begin{lemma}\label{pddimNdim}
Let $M\in A\mathrm{-mod}$. Then $\codim_{\mathcal{F}(N)}(M)$ is
defined and finite if and only if $\mathrm{p.d.}(M)<\infty$.  
\end{lemma}

\begin{proof}
The ``if'' part follows from Proposition~\ref{coresofM} and the fact 
that $H\in\mathcal{F}(N)$.

To prove the ``only if'' part, we take a finite
$\mathcal{F}(N)$-coresolution,
\begin{equation}\label{eqM}
0 \to M \to X_0 \to X_1 \to \dots
\to  X_l \to 0,
\end{equation}
of $M$. From Lemma~\ref{projNN} it follows that all $X_i$ have finite 
projective dimension, which implies that $M$ has finite projective
dimension as well.
\end{proof}

\begin{theorem}\label{Ndeqnablad}
Let $M$ be an $A$-module of finite projective dimension. Then 
\begin{displaymath}
\text{codim}_{\mathcal{F}(N)}(M)=
\text{codim}_{\mathcal{F}(\bnabla)}(M).
\end{displaymath}
\end{theorem}

\begin{proof}
Since $\mathcal{F}(N)\subset\mathcal{F}(\bnabla)$ it follows directly 
that $\text{codim}_{\mathcal{F}(\bnabla)}(M) \leq
\text{codim}_{\mathcal{F}(N)}(M)$.
To prove that $\text{codim}_{\mathcal{F}(N)}(M) \leq
\text{codim}_{\mathcal{F}(\bnabla)}(M)$, we let \eqref{eqM} to be
an $\mathcal{F}(N)$-coresolution of $M$ of minimal length and
$\mathcal{N}_M^{\bullet}$ be the corresponding complex, whose only 
non-zero homology is in degree zero and equals $M$. Applying
Lemma~\ref{apl1} to $\mathtt{A}=A^{\mathrm{opp}}$, 
$X^{(\mathtt{A})}=\DD(\oplus_{i=0}^l X_l)$,
$Y^{(\mathtt{A})}=\DD(T)$ and the complex $\DD(\mathcal{N}_M^{\bullet}[l])$ gives
a complex, the dual $\mathcal{Z}_M^{\bullet}$  of which belongs to 
$\Ka^b(\add(T))$ and is quasi-isomorphic to $\mathcal{N}_M^{\bullet}[l]$.

From a tilting coresolution of $\Delta$ we get the complex  
\begin{displaymath}
\mathcal{J}_{\Delta}^{\bullet}:\quad\quad \dots 
\to 0 \to Q_0 \to Q_1 \to \dots \to Q_r \to 0 \dots
\end{displaymath}
in $\Ka^b(\add(T))$, in which  $Q_0$ is isomorphic to the characteristic
tilting  module $T$, and whose only non-zero homology is in degree
zero and equals $\Delta$.

Let $T(\lambda)$ be a direct summand of $\mathcal{Z}_M^{0}\neq 0$ and $f:T
\twoheadrightarrow T(\lambda) \hookrightarrow \mathcal{Z}_M^{0}$ 
be the projection on this direct summand. By
\cite[Lemma~1]{MO} the homomorphism $f$ gives rise to a non-zero
homomorphism in $\Der^b(A)$ from $\mathcal{J}_{\Delta}^{\bullet}$ to 
$\mathcal{Z}_{M}^{\bullet}$. The last implies that 
$\Ext_A^l(\Delta,M) \neq 0$ and from \cite[Lemma~1]{MP} we obtain that 
$\text{codim}_{\mathcal{F}(\bnabla)}(M)\geq l$. This completes the proof.
\end{proof}

\section{Examples}\label{s8} 

\subsection{Quasi-hereditary algebras}\label{s8.1} 

Assume that $A$ is quasi-hereditary (or, more generally, a properly
stratified algebra, for which the characteristic tilting and cotilting
modules coincide). Then the tilting $A$-module $T$ is also cotilting and 
hence there is an epimorphism,  $T\twoheadrightarrow \nabla$, whose kernel 
is filtered by costandard modules. In this situation Lemma~\ref{propSL} 
implies that $N=\nabla$ and hence $H=I$. Thus, in this case we get that 
$B(A)$ is Morita equivalent to $A$. Moreover, the two-step duality functor 
$G$ is isomorphic to the identity functor. 

\subsection{An algebra with a non-trivial two-step
duality}\label{s8.2}

This example presents an SSS-algebra, $(A,\leq)$, such that the Ringel
dual $(R,\leq_R)$ is a properly stratified algebra with 
$T^{(R)}\neq C^{(R)}$. Let $A$ be the quotient of the path algebra 
of the following quiver
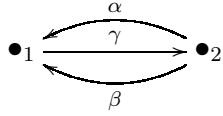
\begin{displaymath}
\xymatrix{
\bullet_{1} \ar[rr]^{\gamma} 
& &  
\bullet_{2}
\ar@/_1pc/[ll]_{\alpha} \ar@/^1pc/[ll]^{\beta}
}
\end{displaymath}
modulo the relations $\gamma\beta=\alpha\gamma\alpha=0$. We set 
$\Lambda=\{1<2\}$. 

The radical filtrations of the projective module $P(\lambda)$, the 
standard module $\Delta(\lambda)$ and the proper standard module 
$\bDelta(\lambda)$, $\lambda=1,2$, look  as follows:
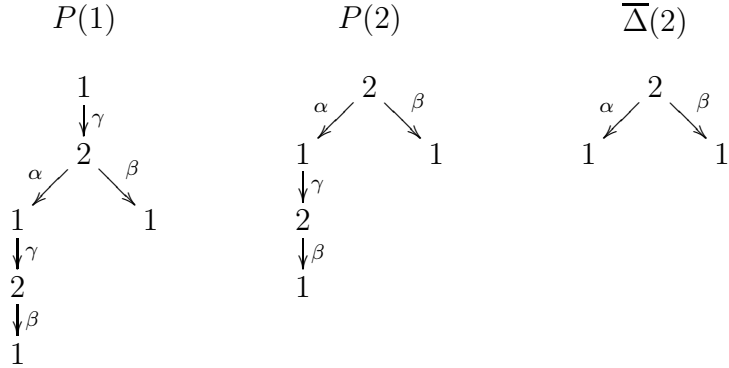
\begin{displaymath}
\xymatrix@!=0.1pc{
  & P(1)  &  \\
  & 1\ar@{->}[d]^{\gamma} &   &  \\
  & 2\ar@{->}[dl]_{\alpha}\ar@{->}[dr]^{\beta} &   &  \\
1\ar@{->}[d]^{\gamma} &   & 1 &  \\
2\ar@{->}[d]^{\beta}&   &   &  \\
1 &   &   &
}
\qquad
\xymatrix@!=0.1pc{
  & P(2)  &  \\
  & 2 \ar@{->}[dl]_{\alpha}\ar@{->}[dr]^{\beta}&   &  \\
1\ar@{->}[d]^{\gamma} &   & 1 &  \\
2\ar@{->}[d]^{\beta} &   &   &  \\
1 &   &   &  \\
  &   &   &
}
\qquad
\xymatrix@!=0.1pc{
  & \bDelta(2)  &  \\
  & 2\ar@{->}[dl]_{\alpha}\ar@{->}[dr]^{\beta} &   &  \\
1 &   & 1 &  \\
  &   &   &  \\
  &   &   &  \\
  &   &   &
}
\end{displaymath}
and  $\Delta(2)=P(2)$, $\Delta(1)=\bDelta(1)=L(1)$. Here we see that
$A$ is an SSS-algebra, but not properly stratified.

The injective module $I(\lambda)$, the costandard module
$\nabla(\lambda)$ and the proper costandard module $\bnabla(\lambda)$, 
$\lambda=1,2$, have the following socle filtrations:
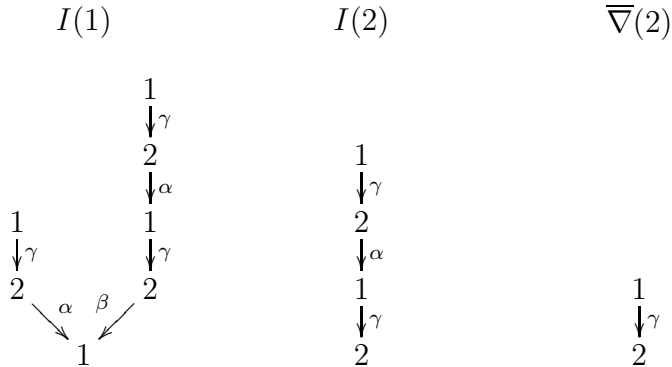
\begin{displaymath}
\xymatrix@!=0.1pc{
  & I(1)  &  \\
  &   & 1\ar@{->}[d]^{\gamma} &  \\
  &   & 2\ar@{->}[d]^{\alpha} &  \\
1\ar@{->}[d]^{\gamma} &   & 1\ar@{->}[d]^{\gamma} &  \\
2\ar@{->}[dr]^{\alpha} &   & 2\ar@{->}[dl]_{\beta} &  \\
  & 1 &   &
}
\qquad
\xymatrix@!=0.1pc{
  & I(2)  &  \\
  &   &   &  \\
  & 1\ar@{->}[d]^{\gamma} &   &  \\
  & 2\ar@{->}[d]^{\alpha} &   &  \\
  & 1\ar@{->}[d]^{\gamma} &   &  \\
  & 2 &   &
}
\qquad
\xymatrix@!=0.1pc{
  & \bnabla(2)  & \\
  &   &   &  \\
  &   &   &  \\
  &   &   &  \\
  & 1 \ar@{->}[d]^{\gamma} &   &  \\
  & 2 &   &
}
\end{displaymath}
and $\nabla(2)=I(2)$, $\nabla(1)=\bnabla(1)=L(1)$. 

The modules $T(\lambda)$, $N(\lambda)$, $S(\lambda)$ and 
$H(\lambda)$ have the following radical filtrations:
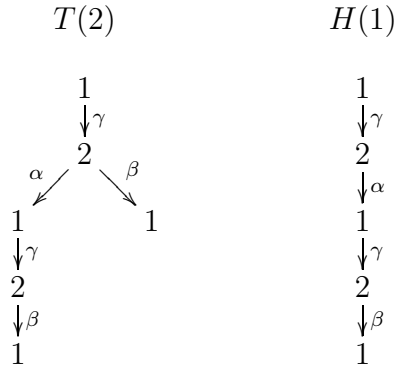
\begin{displaymath}
\xymatrix@!=0.1pc{
  & T(2)  & \\
  & 1\ar@{->}[d]^{\gamma} &   & \\
  & 2\ar@{->}[dl]_{\alpha}\ar@{->}[dr]^{\beta} &   & \\
1\ar@{->}[d]^{\gamma} &   & 1 & \\
2\ar@{->}[d]^{\beta} &   &   & \\
1 &   &   &
}
\qquad
\xymatrix@!=0.1pc{
  & H(1)  & \\
  & 1\ar@{->}[d]^{\gamma} &   & \\
  & 2\ar@{->}[d]^{\alpha} &   & \\
  & 1\ar@{->}[d]^{\gamma} &   & \\
  & 2\ar@{->}[d]^{\beta} &   & \\
  & 1 &   &
}
\end{displaymath}
and $T(1)=N(1)=L(1)$, $N(2)=\nabla(2)$, $S(1)=0$, $S(2)=L(1)\oplus
L(1)$, $H(2)=I(2)$.

Since $S(1)$ and $S(2)$ are in $\mathcal{F}(N)$, we can conclude, by 
Theorem~\ref{Ringelpropstr}, that the Ringel dual $(R,\leq_R)$ is
properly stratified. Also $T^{(R)}\neq C^{(R)}$, because $H(1)\neq
I(1)$. 

By a straightforward calculation one gets that the Ringel dual $R$ is
the quotient of the path algebra of the following quiver
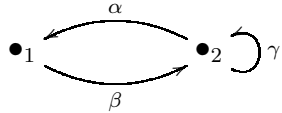
\begin{displaymath}
\xymatrix{
\bullet_{1} \ar@/_1pc/[rr]_{\beta} & &  
\bullet_{2}
\ar@/_1pc/[ll]_{\alpha}\ar@(dr,ur)[]_{\gamma}
} 
\end{displaymath}
modulo the relations $\gamma\beta=\gamma^2=\alpha\beta = 0$. The projective modules
over this algebra have the following radical filtrations:
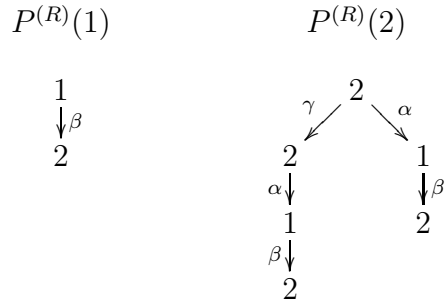
\begin{displaymath}
\xymatrix@!=0.1pc{
  & P^{(R)}(1)  & \\
  & 1\ar@{->}[d]^{\beta} & \\
  & 2& \\
}
\qquad\qquad
\xymatrix@!=0.1pc{
   & P^{(R)}(2)  & \\
   & 2\ar@{->}[dl]_{\gamma}\ar@{->}[dr]^{\alpha}     &   & \\
 2\ar@{->}[d]_{\alpha} &       & 1\ar@{->}[d]^{\beta}  & \\
 1\ar@{->}[d]_{\beta} &       & 2  & \\
 2 &       &    & \\
}
\end{displaymath}
This algebra has dimension $8$ and is properly stratified with 
respect to the opposite order $2<_R 1$ on $\Lambda$.

By a straightforward calculation one gets that the  two-step dual 
algebra $B(A)$ is the quotient of the path algebra of the following 
quiver
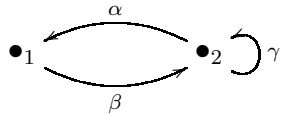
\begin{displaymath}
\xymatrix{
\bullet_{1} \ar@/_1pc/[rr]_{\beta}& &  
\bullet_{2}
 \ar@/_1pc/[ll]_{\alpha} \ar@(dr,ur)[]_{\gamma}
}
\end{displaymath}
modulo the relations $\gamma^2=\gamma\beta=\beta\alpha=0$. The projective modules
over this algebra have the following radical filtrations:
\begin{displaymath}
\xymatrix@!=0.1pc{
  & P^{(B(A))}(1)  & \\
  & 1\ar@{->}[d]^{\beta} & \\
  & 2\ar@{->}[d]^{\alpha}& \\
  & 1& \\
}
\qquad\qquad
\xymatrix@!=0.1pc{
   & P^{(B(A))}(2)  & \\
   & 2\ar@{->}[dl]_{\gamma}\ar@{->}[dr]^{\alpha}     &   & \\
 2\ar@{->}[d]_{\alpha} &       & 1  & \\
 1 &       &   & \\
}
\end{displaymath}
The dimension of $B(A)$ is $7$, while the dimension of $A$ is $11$. 
Since both algebras are basic, it follows that $A$ and $B(A)^{\mathrm{opp}}$ 
are neither isomorphic nor Morita equivalent.

\subsection{A properly stratified algebra with duality whose two-step
dual is not properly stratified}\label{s8.25}

This example presents a properly stratified algebra $(A,\leq)$ 
having a simple preserving duality, ${}^\circ$, whose  
Ringel dual $(R,\leq_R)$ is properly stratified, does not have
any simple preserving duality and such that the two-step dual algebra 
$B(A)$ is not properly stratified.

Let $A$ be the quotient of the path algebra of the following quiver
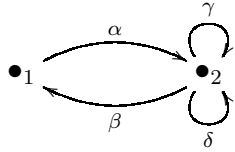
\begin{displaymath}
\xymatrix{
\bullet_{1} \ar@/^1pc/[rr]^{\alpha} 
& &  
\bullet_{2}
\ar@/^1pc/[ll]^{\beta} \ar@(ul,ur)[]^{\gamma}  \ar@(dl,dr)[]_{\delta} 
}
\end{displaymath}
modulo the relations $\gamma^2=\delta^2=\gamma\delta=\delta\gamma=
\alpha\beta=0$. We set $\Lambda=\{1<2\}$. The algebra $A$ is
isomorphic to the opposite algebra via the antiautomorphism 
$\iota:A\to A$, defined by $\iota(e_i)=e_i$, $i=1,2$,
$\iota(\alpha)=\beta$, $\iota(\beta)=\alpha$, $\iota(\gamma)=\delta$, 
$\iota(\delta)=\gamma$. Since $\iota$ stabilizes the primitive
idempotents, it induces a simple preserving duality for $A$.

The radical filtrations of the projective module $P(\lambda)$, the 
standard module $\Delta(\lambda)$ and the proper standard module 
$\bDelta(\lambda)$, $\lambda=1,2$, look as follows:
\begin{displaymath}
\xymatrix@!=0.1pc{
  & P(1)  &  \\
  & 1\ar@{->}[d]^{\alpha} &   &  \\
  & 2\ar@{->}[d]^{\beta}
\ar@{->}[dl]_{\gamma}\ar@{->}[dr]^{\delta} &   &  \\
2\ar@{->}[d]^{\beta} & 1 & 2\ar@{->}[d]^{\beta} &  \\
1 &   & 1 &  
}
\qquad
\xymatrix@!=0.1pc{
  & P(2)  &  \\
  & 2 \ar@{->}[d]^{\beta} 
\ar@{->}[dl]_{\gamma} \ar@{->}[dr]^{\delta}  &  \\
2\ar@{->}[d]^{\beta} & 1 & 2\ar@{->}[d]^{\beta} &  \\
1 &   & 1 &  \\
  &   &   &  
}
\qquad
\xymatrix@!=0.1pc{
  & \bDelta(2)  &  \\
  & 2\ar@{->}[d]^{\beta} &   &  \\
  & 1 &   &  \\
  &   &   &  \\
  &   &   &  \\
  &   &   &
}
\end{displaymath}
and  $\Delta(2)=P(2)$, $\Delta(1)=\bDelta(1)=L(1)$. It follows that 
$A$ is properly stratified.

The injective module $I(\lambda)$, the costandard module
$\nabla(\lambda)$ and the proper costandard module 
$\bnabla(\lambda)$, $\lambda=1,2$, have the following socle 
filtrations (dual to the corresponding radical filtrations above):
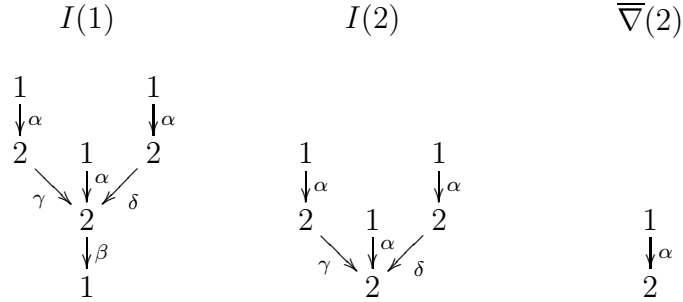
\begin{displaymath}
\xymatrix@!=0.1pc{
  & I(1)  &  \\
1\ar@{->}[d]^{\alpha} &   & 1 \ar@{->}[d]^{\alpha}& \\
2\ar@{->}[dr]_{\gamma} & 1\ar@{->}[d]^>>>>{\alpha} & 2 
\ar@{->}[dl]^{\delta}& \\
  & 2\ar@{->}[d]^{\beta} &   & \\
  & 1 &   &
}
\qquad
\xymatrix@!=0.1pc{
  & I(2)  &  \\
  &   &   &  \\
1 \ar@{->}[d]^{\alpha}&   & 1 \ar@{->}[d]^{\alpha}& \\
2\ar@{->}[dr]_{\gamma} & 1\ar@{->}[d]^>>>>{\alpha} & 2 
\ar@{->}[dl]^{\delta}& \\
  & 2 &   &
}
\qquad
\xymatrix@!=0.1pc{
  & \bnabla(2)  & \\
  &   &   &  \\
  &   &   &  \\
  & 1 \ar@{->}[d]^{\alpha} &   &  \\
  & 2 &   &
}
\end{displaymath}
and $\nabla(2)=I(2)$, $\nabla(1)=\bnabla(1)=L(1)$. 

The tilting module $T(2)$ has the following radical filtration:
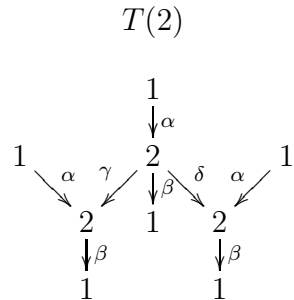
\begin{displaymath}
\xymatrix@!=0.1pc{
 &   & T(2)  &  \\
 &   &  1\ar@{->}[d]^{\alpha}&   &  \\
1\ar@{->}[dr]^{\alpha} &  &
 2\ar@{->}[dl]_{\gamma}\ar@{->}[dr]^{\delta}
\ar@{->}[d]^{\beta} &
   & 1 \ar@{->}[dl]_{\alpha}\\
 & 2\ar@{->}[d]^{\beta} & 1  & 
2\ar@{->}[d]^{\beta} &  \\
 & 1 &   & 1 &  \\
}
\end{displaymath}
and $T(1)=L(1)$.

By a straightforward calculation one gets that the Ringel dual $R$ is
the  quotient of the path algebra of the following quiver
\begin{displaymath}
\xymatrix{
\bullet_{2}
\ar[rr]_{c}\ar@(ur,ul)[]_{d}
\ar@(dr,dl)[]^{e} & &  
\bullet_{1}\ar@/_1pc/[ll]_{a}\ar@/^1pc/[ll]^{b}
} 
\end{displaymath}
modulo the relations $d^2=e^2=d e=e d= e a=d b=c a=c b =c d a = c e b =0$
and $d a=e b$. The projective modules
over this algebra have the following radical filtrations:
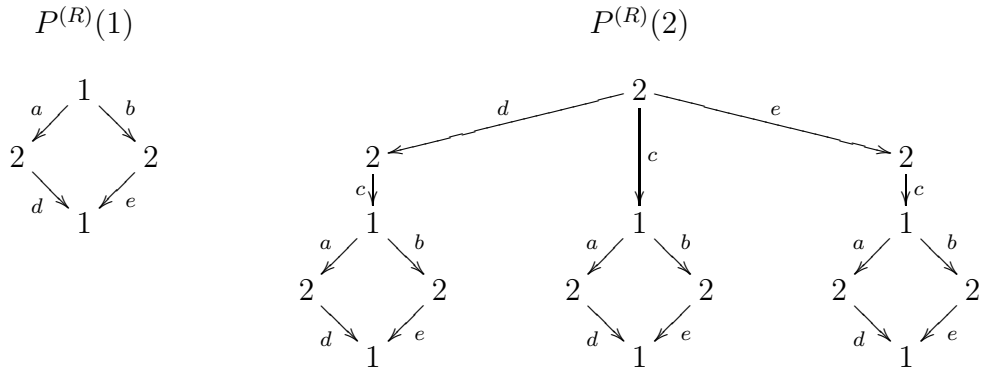
\begin{displaymath}
\xymatrix@!=0.1pc{
  & P^{(R)}(1)  & \\
  & 1\ar@{->}[dl]_{a}\ar@{->}[dr]^{b} & \\
 2\ar@{->}[dr]_{d} & & 2\ar@{->}[dl]^{e} \\
  & 1& \\
}
\qquad\qquad
\xymatrix@!=0.1pc{
&&&&& P^{(R)}(2)  &&&&& \\
&&&&&2\ar@{->}[dllll]_{d}\ar@{->}[drrrr]^{e}\ar@{->}[dd]^{c}&&&&& \\
&2\ar@{->}[d]_{c}&&&&&&&&2\ar@{->}[d]^{c}& \\
&1\ar@{->}[dl]_{a}\ar@{->}[dr]^{b}&&&&1\ar@{->}[dl]_{a}\ar@{->}[dr]^{b}
&&&&1\ar@{->}[dl]_{a}\ar@{->}[dr]^{b}& \\
2\ar@{->}[dr]_{d}&&2\ar@{->}[dl]^{e}&&2\ar@{->}[dr]_{d}&&2\ar@{->}[dl]^{e}&&
2\ar@{->}[dr]_{d}&&2\ar@{->}[dl]^{e} \\
&1&&&&1&&&&1& \\
}
\end{displaymath}
This algebra is properly stratified with respect to the opposite order 
$2<_R 1$ on $\Lambda$. Further it is easy to get the following
equalities: $\text{dim}_{\k}\Ext_R^1(L^{(R)}(1),L^{(R)}(2))=2$ and   
$\text{dim}_{\k}\Ext_R^1(L^{(R)}(2),L^{(R)}(1))=1$ and hence $R$ does
not have any simple preserving duality. 

By a straightforward calculation one also gets that the two-step dual 
algebra $B(A)$ is the quotient of the path algebra of the following 
quiver
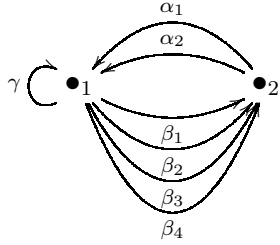
\begin{displaymath}
\xymatrix{
\bullet_{1} \ar@/_1pc/[rr]_{\beta_1}\ar@/_2pc/[rr]_{\beta_2}
\ar@/_3pc/[rr]_{\beta_3}\ar@/_4pc/[rr]_{\beta_4} 
\ar@(dl,ul)[]^{\gamma} & &  
\bullet_{2}
 \ar@/_1pc/[ll]_{\alpha_2} \ar@/_2pc/[ll]_{\alpha_1} 
}
\end{displaymath}
modulo the relations
$\beta_1\gamma=\beta_2\gamma=\beta_3\gamma=\beta_4\gamma=
\gamma\alpha_1=\gamma\alpha_2=0$,
$\beta_3\alpha_1\beta_2=\beta_4\alpha_2\beta_1=0$,
$\beta_3\alpha_1\beta_1=\beta_4\alpha_2\beta_2$,
$\beta_3\alpha_1\beta_4=\beta_4\alpha_2\beta_4=
\beta_3\alpha_1\beta_3=\beta_4\alpha_2\beta_3=0$,
$\beta_1\alpha_1=\beta_2\alpha_1=\beta_1\alpha_2=
\beta_2\alpha_2=\beta_3\alpha_2=\beta_4\alpha_1=0$ and $\gamma^2=0$. 
The projective modules over this algebra have the following radical 
filtrations:
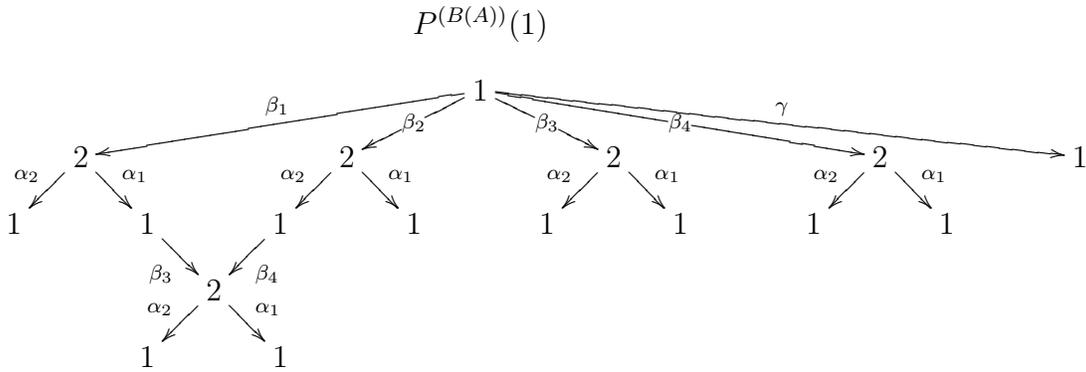
\begin{gather*}
\xymatrix@!=0.1pc{
&&&&&&& P^{(B(A))}(1)  &&&&&&&&& \\
&&&&&&& 1\ar[dllllll]_{\beta_1}\ar[dll]|-{\beta_2}\ar[drr]|-{\beta_3}
\ar[drrrrrr]|-{\beta_4}\ar[drrrrrrrrr]^{\gamma} &&&&&&&&&\\
&2\ar[dl]_{\alpha_2}\ar[dr]^{\alpha_1}&&&&2\ar[dl]_{\alpha_2}\ar[dr]^{\alpha_1}
&&&&2\ar[dl]_{\alpha_2}\ar[dr]^{\alpha_1}&&&&2\ar[dl]_{\alpha_2}\ar[dr]^{\alpha_1}&&&1\\
1&&1\ar[dr]_{\beta_3}&&1\ar[dl]^{\beta_4}&&1&&1&&1&&1&&1&&\\
&&&2\ar[dl]_{\alpha_2}\ar[dr]^{\alpha_1}&&&&&&&&&&&&&\\
&&1&&1&&&&&&&&&&&&\\
}\\
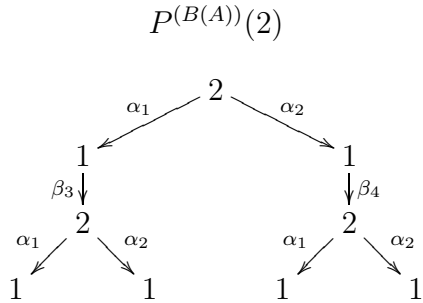
\xymatrix@!=0.1pc{
&&&P^{(B(A))}(2)&&&\\
&&&2\ar[dll]_{\alpha_1}\ar[drr]^{\alpha_2}&&&\\
&1\ar[d]_{\beta_3}&&&&1\ar[d]^{\beta_4}&\\
&2\ar[dl]_{\alpha_1}\ar[dr]^{\alpha_2}&&&&2\ar[dl]_{\alpha_1}\ar[dr]^{\alpha_2}&\\
1&&1&&1&&1\\
}
\end{gather*}
One immediately sees that $B(A)^\mathrm{opp}$ is an SSS-algebra, while 
$B(A)$ is not (since the trace of $P^{(B(A))}(2)$ in $P^{(B(A))}(1)$ is not a direct
sum of several copies of $P^{(B(A))}(2)$). In particular, $B(A)$ is neither properly 
stratified nor has a simple preserving duality.

\subsection{Stratified algebras associated with Harish-Chandra
bimodules}\label{s8.3}

Let $\mathfrak{g}$ be a semi-simple finite dimensional complex Lie
algebra, $U(\mathfrak{g})$ be its universal enveloping algebra and  
$Z(\mathfrak{g})$ be the center of $U(\mathfrak{g})$. Denote by
$\mathcal{H}$
the category of all Harish-Chandra $\mathfrak{g}$-bimodules, that is
finitely generated $\mathfrak{g}$-bimodules, which are direct sums of
finite dimensional $\mathfrak{g}$-modules under the adjoint action of 
$\mathfrak{g}$. Fix a positive integer $n$ and two maximal ideals
$\chi$ and $\xi$ of $Z(\mathfrak{g})$. Suppose for simplicity that 
$\chi$ and $\xi$ correspond to regular central characters and
denote by ${}_\chi \mathcal{H}_\xi^n$ the full subcategory in
$\mathcal{H}$ which consists of all bimodules $M$ satisfying 
$M \xi^n = 0$ and $\chi^m M =0$, $m\gg 0$. We refer the reader to
\cite[Kapitel~6]{Ja} and  \cite{So} for details.  
It is well-known, see for example \cite[Section~5]{So}, that 
${}_\chi \mathcal{H}_\xi^n$ is equivalent to the module category of
some finite dimensional associative algebra $A$. In
\cite[Section~5.6]{Ma} it was shown that $A$ is properly stratified. 
At the same time, if $\rank(\mathfrak{g})>1$ and $n>1$, it is easy 
to see that the cotilting modules for $A$ have infinite projective 
dimension and hence are not tilting modules. 

However, using the translation functors one can easily show that the
Ringel dual of $A$ is properly stratified. The simple bimodules in 
${}_\chi \mathcal{H}_\xi^n$ are indexed by the elements of the Weyl
group $W$, see \cite[Kapitel~6]{Ja}. Let $w_0$ be the longest element 
in $W$. Consider the tilting bimodule $T(w_0)$ corresponding to
$w_0$.  Then all tilting bimodules in ${}_\chi \mathcal{H}_\xi^n$ are 
direct summands of $F\otimes T(w_0)$ for some finite dimensional 
$\mathfrak{g}$-module $F$. Let us construct the bimodules $N(w)$, 
$w\in W$, inductively. Set $N(w_0)=T(w_0)$, and for $w\in W$ and a 
simple reflection $s$ such that the length of $sw$ is smaller than 
the length of $w$, let $N(sw)$ be the cokernel of the adjunction 
morphism from $N(w)$ to $\theta_s(N(w))$, where $\theta_s$ denotes 
the translation functor through the $s$-wall. It is easy to see  
that this adjunction morphism is always injective in the situation 
above. Using the standard properties of the translation functors and 
the fact that all tilting bimodules are obtained by translating 
$T(w_0)$, we obtain that all tilting bimodules are filtered by $N(w)$, 
$w\in W$. This implies that the Ringel dual $R$ of $A$ is properly 
stratified. 
     
In fact, using Arkhipov's functor, \cite{Ar}, one can show that $R$ is
isomorphic to $A$. Iterating this we also obtain $A\cong B(A)$.
However, in this case the two-step duality functor is quite far from
being trivial, because the injective $A$-modules have infinite
projective dimension in general. In fact, the two-step duality
functor does something remarkable, namely, it defines a covariant
equivalence between the categories of $\mathcal{P}(A)^{<\infty}$ and  
$\mathcal{I}(A)^{<\infty}$. That these two categories are
contravariantly equivalent follows from the existence of a simple
preserving duality for $A$ (the last comes from the equivalence of 
${}_\chi \mathcal{H}_\xi^n$ with a block of the thick category
$\mathcal{O}$, \cite[Theorem~3]{So}). In particular, the category 
$\mathcal{P}(A)^{<\infty}$ happens to be equivalent to  
$\big(\mathcal{P}(A)^{<\infty}\big)^\mathrm{opp}$. Obviously, this
gives us a contravariant equivalence from $\mathcal{P}(A)^{<\infty}$   
to itself, which sends $P$ to $H$. In particular, from
Proposition~\ref{c63.25} it follows that the finitistic dimension of 
${}_\chi \mathcal{H}_\xi^n$ equals twice the projective dimension of
the characteristic tilting module in ${}_\chi \mathcal{H}_\xi^n$.

\subsection{A tensor construction for properly stratified algebras}
\label{s8.4}

In this section we present one general construction of properly
stratified algebras using quasi-hereditary and local algebras. In 
this way one obtains a huge family of stratified algebras, moreover, 
there is an easy criterion when all tilting modules over such algebras 
are cotilting. In fact, this happens to be the case if and only if the 
local algebra we start from is self-injective. Thus this gives us a 
possibility of constructing series of properly stratified algebras for 
which tilting and cotilting  modules do not coincide, hence showing 
that the two-step duality we worked out is not that rare.   

We remark that any local algebra is properly stratified with simple
proper standard and costandard modules, projective standard modules
and injective costandard modules.

Throughout this subsection we fix a quasi-hereditary algebra, 
$(\mathtt{A},\leq)$, with $\Lambda$ being an index set for the isomorphism 
classes of simple $\mathtt{A}$-modules. Let further $\mathtt{B}$ be a fixed 
local algebra and consider the algebra $D = \mathtt{A}\otimes_{\k} \mathtt{B}$. 
Since $\mathtt{B}$ is local, $\Lambda$ also indexes the isomorphism classes 
of simple $D$-modules in a natural way. Note that we have the following:
$\Delta^{(\mathtt{A})}(\lambda)=\bDelta^{(\mathtt{A})}(\lambda)$,
$\nabla^{(\mathtt{A})}(\lambda)=\bnabla^{(\mathtt{A})}(\lambda)$ for all $\lambda$.
Moreover,  $P^{(\mathtt{B})}=\Delta^{(\mathtt{B})}$, 
$I^{(\mathtt{B})}=\nabla^{(\mathtt{B})}$, 
$\bDelta^{(\mathtt{B})}=\bnabla^{(\mathtt{B})}=
L^{(\mathtt{B})}$.

\begin{proposition}\label{TCT0}
The algebra $(D,\leq)$ is properly stratified. Moreover, for each 
$\lambda \in \Lambda$ we have 
\begin{displaymath}
X^{(D)}(\lambda) =  X^{(\mathtt{A})}(\lambda){\otimes_{\k}} X^{(\mathtt{B})},
\quad\text{where}\quad 
X\in\{\Delta,\bDelta,\nabla,\bnabla\}.
\end{displaymath}
\end{proposition}

\begin{proof}
Let $\lambda\in \Lambda$. Then we obtain 
\begin{displaymath}
P^{(D)}(\lambda)= D(e_\lambda{\otimes_{\k}} 1_\mathtt{B}) = 
\mathtt{A}e_\lambda {\otimes_{\k}} \mathtt{B} 1_\mathtt{B} =
P(\lambda){\otimes_{\k}} P^{(\mathtt{B})},  
\end{displaymath}
where the idempotent $e_\lambda\in \mathtt{A}$ is chosen such that
$P^{(\mathtt{A})}(\lambda)=\mathtt{A} e_\lambda$. The functor 
${}_-{\otimes_{\k}} P^{(\mathtt{B})}:\mathtt{A}\text{-mod}\to D\text{-mod}$ is exact 
(as the tensor product over a field) and if we apply 
${}_-{\otimes_{\k}} P^{(\mathtt{B})}$ to the short exact sequence 
\begin{displaymath}
0 \to K^{(\mathtt{A})} \to P^{(\mathtt{A})}(\lambda) \to 
\Delta^{(\mathtt{A})}(\lambda) \to 0,
\end{displaymath}
where $K^{(\mathtt{A})}$ has a filtration with subquotients 
$\Delta^{(\mathtt{A})}(\mu)$, $\mu>\lambda$, we obtain the short exact 
sequence
\begin{displaymath}
0 \to K^{(\mathtt{A})}{\otimes_{\k}} P^{(\mathtt{B})} \to 
P^{(D)}(\lambda) \to 
\Delta^{(\mathtt{A})}(\lambda){\otimes_{\k}} P^{(\mathtt{B})} \to 0,
\end{displaymath}
where, using the exactness of ${}_-{\otimes_{\k}} P^{(\mathtt{B})}$, 
$K^{(\mathtt{A})}{\otimes_{\k}} P^{(\mathtt{B})}$ has a filtration with subquotients  
$\Delta^{(\mathtt{A})}(\mu){\otimes_{\k}} P^{(\mathtt{B})}$, $\mu>\lambda$. 
Exactness of ${}_-{\otimes_{\k}} P^{(\mathtt{B})}$, $L^{(\mathtt{A})}(\nu)
{\otimes_{\k}} {}_-$, $\nu\in\Lambda$, and the equality 
$L^{(D)}(\mu)=L^{(\mathtt{A})}(\mu)\otimes_{\k}L^{(\mathtt{B})}$ also implies that 
$[\Delta^{(\mathtt{A})}(\lambda){\otimes_{\k}} P^{(\mathtt{B})}:L^{(D)}(\mu)]=0$ 
for $\mu>\lambda$.  Therefore we conclude that 
the property (SS) holds for $\mathtt{A}\otimes_{\k} \mathtt{B}$. 

To show that (PS) holds, apply the exact functor 
$\Delta^{(\mathtt{A})}(\lambda){\otimes_{\k}} {}_-:\mathtt{B}\text{-mod}\to 
D\text{-mod}$ to the short exact sequence 
\begin{displaymath}
0 \to K^{(\mathtt{B})} \to P^{(\mathtt{B})} \to L^{(\mathtt{B})} \to 0,
\end{displaymath}
where $K^{(\mathtt{B})}$ has a filtration with simple subquotients 
$L^{(\mathtt{B})}$, and obtain  
\begin{displaymath}
0 \to \Delta^{(\mathtt{A})}(\lambda){\otimes_{\k}} K^{(\mathtt{B})} \to
\Delta^{(\mathtt{A})}(\lambda){\otimes_{\k}} P^{(\mathtt{B})} \to 
\Delta^{(\mathtt{A})}(\lambda){\otimes_{\k}} L^{(\mathtt{B})} \to 0,
\end{displaymath}
where $\Delta^{(\mathtt{A})}(\lambda)\otimes K^{(\mathtt{B})}$ has a filtration 
with subquotients $\Delta^{(\mathtt{A})}(\lambda){\otimes_{\k}} L^{(\mathtt{B})}$. 
We conclude that property (PS) holds. By dual arguments it follows that 
$I^{(D)}(\lambda)= I^{(\mathtt{A})}(\lambda){\otimes_{\k}} I^{(\mathtt{B})}$ 
has a filtration with subquotients $\nabla^{(D)} =  \nabla^{(\mathtt{A})} 
{\otimes_{\k}} \nabla^{(\mathtt{B})}$, and $\nabla^{(D)}$ has a filtration with 
subquotients $\bnabla^{(D)}= \nabla^{(\mathtt{A})} {\otimes_{\k}} L^{(\mathtt{B})}$. 
Hence the lemma is proved.
\end{proof}

As an immediate corollary we obtain:
 
\begin{corollary}\label{TCT1}
The module $\Delta^{(D)}$ is projective and the module $\nabla^{(D)}$ 
is injective if viewed as $\mathtt{B}$-module.
\end{corollary}

\begin{corollary}\label{TCT15}
We have:
\begin{enumerate}
\item $T^{(D)}(\lambda) =  T^{(\mathtt{A})}(\lambda)
 {\otimes_{\k}} T^{(\mathtt{B})} =  T^{(\mathtt{A})}(\lambda)
 {\otimes_{\k}} P^{(\mathtt{B})}$;
\item $C^{(D)}(\lambda) =  C^{(\mathtt{A})}(\lambda)
 {\otimes_{\k}} C^{(\mathtt{B})} =  T^{(\mathtt{A})}(\lambda)
 {\otimes_{\k}} I^{(\mathtt{B})} $.
\end{enumerate}
\end{corollary}

\begin{proof}
Follows from Proposition~\ref{TCT0} and exactness of the functors
${}_-{\otimes_{\k}} P^{(\mathtt{B})}$ and ${}_-{\otimes_{\k}} I^{(\mathtt{B})}$.
\end{proof}

And, finally, we can formulate our main result in this section:

\begin{corollary}\label{TCT2}
The tilting modules for $D$ are cotilting if and only if $\mathtt{B}$ is 
self-injective.
\end{corollary}

\begin{proof}
Since $D$ is properly stratified, we have that $T^{(D)}=C^{(D)}$
is equivalent to $C^{(D)}\in\mathcal{F}(\Delta^{(D)})$. By 
Lemma~\ref{identdelta} the last is equivalent to 
$\text{p.d.}(C^{(D)})<\infty$. From Corollary~\ref{TCT1} and the 
fact that $D$ is properly stratified we get that all projective 
$D$-modules are also projective as $\mathtt{B}$-modules.

Thus  $\text{p.d.}(C^{(D)})<\infty$ implies that $C^{(D)}$ is a
$\mathtt{B}$-module of finite projective dimension, hence projective as 
$\mathtt{B}$ is local. But from Corollary~\ref{TCT1} it also follows that 
$C^{(D)}$ is $\mathtt{B}$  injective and hence $\mathtt{B}$ is self-injective. 

On the other hand, Corollary~\ref{TCT15} implies that, when $\mathtt{B}$ is 
self-injective, then
\begin{displaymath}
T^{(D)}=T^{(\mathtt{A})}\otimes_{\k} T^{(\mathtt{B})} = 
T^{(\mathtt{A})}\otimes_{\k} P^{(\mathtt{B})} = 
C^{(\mathtt{A})}\otimes_{\k} I^{(\mathtt{B})} = 
C^{(\mathtt{A})}\otimes_{\k} C^{(\mathtt{B})} =C^{(D)}.
\end{displaymath}
This completes the proof.
\end{proof}

\section{Appendix: Two technical lemmas}\label{sapp}

In this appendix we prove two auxiliary technical lemmas similar 
to \cite[Lemma~4]{MO}, which were used in the paper.

\begin{lemma}\label{apl1}
Let $\mathtt{A}$ be a finite-dimensional associative $\Bbbk$-algebra,
$X^{(\mathtt{A})}$ be a $\mathtt{A}$-module, and $Y^{(\mathtt{A})}$ be 
a (generalized) (co)tilting $\mathtt{A}$-module. 
Assume that $X^{(\mathtt{A})}$ has a finite coresolution by modules from 
$\mathrm{add}(Y^{(\mathtt{A})})$. Let $\mathcal{X}^{\bullet}$ be a positive 
complex in $\Ka(\add(X^{(\mathtt{A})}))$. Then there is a positive complex,
$\mathcal{Y}^{\bullet}$, in $\Ka(\add(Y^{(\mathtt{A})}))$ such that
$\mathcal{X}^{\bullet}$ is quasi-isomorphic to
$\mathcal{Y}^{\bullet}$. Moreover, if $\mathcal{X}^{\bullet} \in
\Ka^{b}(\add(X^{(\mathtt{A})}))$, then $\mathcal{Y}^{\bullet} \in 
\Ka^{b}(\add(Y^{(\mathtt{A})}))$. 
\end{lemma}

\begin{proof}
Let  
\begin{displaymath}
\mathcal{X}^{\bullet}: \quad\quad \dots \to 0 \rightarrow 
X^{(\mathtt{A})}_{0} \rightarrow X^{(\mathtt{A})}_{1} \rightarrow \dots
\end{displaymath}
be a positive complex in $\Ka(\add(X^{(\mathtt{A})}))$ and put
$\mathcal{X}_{j}^{\bullet} = t_{j}\mathcal{X}^{\bullet}$. We will show 
by induction that for each $j \geq 0$ there exist a complex,
$\mathcal{Y}_{j}^{\bullet} \in \Ka(\add(Y^{(\mathtt{A})}))$, and a
quasi-isomorphism, $\Phi_{j}:\mathcal{X}_{j}^{\bullet} \rightarrow
\mathcal{Y}_{j}^{\bullet}$. Moreover, we will choose the family
$\{\mathcal{Y}_{j}^{\bullet} \}_{j\geq 0}$ such that for all $ k < j$
we have $t_{k}\mathcal{Y}_{j}^{\bullet} =
t_k\mathcal{Y}_{k}^{\bullet}$. 

Since $X^{(\mathtt{A})}_{i}\in \add(X^{(\mathtt{A})})$ for all $i$ we 
can choose by our assumptions a coresolution, $\mathcal{Z}_{i}^{\bullet} 
\in \Ka(\add(Y^{(\mathtt{A})}))$, of $X^{(\mathtt{A})}_{i}$, and a quasi-isomorphism, 
$\phi_{i}:(X^{(\mathtt{A})}_{i})^{\bullet} \rightarrow \mathcal{Y}_{i}^{\bullet}$. 
In the case $j = 0$ we put $\mathcal{Y}_{0}^{\bullet}=\mathcal{Z}_{0}^{\bullet}$ 
and we are done.   

Now suppose by induction that there exists a quasi-isomorphism,
$\Phi_{j-1}:\mathcal{X}_{j-1}^{\bullet} \rightarrow 
\mathcal{Y}_{j-1}^{\bullet}$. The map $d_{j-1}:X^{(\mathtt{A})}_{j-1}\rightarrow 
X^{(\mathtt{A})}_{j}$ induces the distinguished triangle 
\begin{displaymath}
\mathcal{X}_{j-1}^{\bullet} 
\overset{d_{j-1}^{\bullet}}{\longrightarrow} (X^{(\mathtt{A})}_{j})^{\bullet}[-j+1] 
\rightarrow \text{Cone}(d_{j-1}^{\bullet}) \rightarrow 
\mathcal{X}_{j-1}^{\bullet}[1],
\end{displaymath}
and we have $\text{Cone}(d_{j-1}^{\bullet})
=\mathcal{X}_{j}^{\bullet}[1]$. 

Using $\mathrm{Ext}_{\mathtt{A}}(Y^{(\mathtt{A})},Y^{(\mathtt{A})})=0$ and 
\cite[Chap.III, Lemma~2.1]{Ha} we can choose a representative, 
$\psi_{j}:\mathcal{Y}_{j-1}^{\bullet}
\rightarrow \mathcal{Z}_{j}^{\bullet}[-j+1]$, in
$\Ka(\add(Y^{(\mathtt{A})}))$ of the composition
$\phi_{j}[-j+1] \circ  d_{j-1}^\bullet \circ \Phi_{j-1}^{-1}$ (in
$\Der^{b}(\mathtt{A}\text{-mod})$). This gives us a diagram in 
$\Der^{b}(\mathtt{A}\text{-mod})$, which can be completed to the following  
commutative diagram:  
\begin{displaymath}
\xymatrix{
\mathcal{X}_{j-1}^{\bullet} \ar[r]^>>>>>{d_{j-1}^{\bullet}}
\ar[d]^{\Phi_{j-1}}_{\wr} &
(X^{(\mathtt{A})}_{j})^{\bullet}[-j+1] \ar[r]\ar[d]^{\phi_{j}[-j+1]}_{\wr} &
\mathcal{X}_{j}^{\bullet}[1] \ar[r]\ar[d]^{\Phi_{j}} &
\mathcal{X}_{j-1}^{\bullet}[1] \ar[d]^{\Phi_{j-1}[1]}\\ 
\mathcal{Y}_{j-1}^{\bullet} \ar[r]^>>>>>{\psi_{j}} & 
\mathcal{Z}_{j}^{\bullet}[-j+1]
\ar[r] & \text{Cone}(\psi_{j}) \ar[r] & 
\mathcal{Y}_{j-1}^{\bullet}[1]. \\
}
\end{displaymath}
Since both $\Phi_{j-1}$ and  $\phi_{j}[-j+1]$ are isomorphisms 
the morphism $\Phi_{j}$ is an isomorphism too. Hence we get the 
quasi-isomorphism $\Phi_{j}:\mathcal{X}_{j}^{\bullet} \to 
\text{Cone}(\psi_{j})[-1]$, where $\text{Cone}(\psi_{j})[-1]$ is a
positive complex in $\Ka^{b}(\add(Y^{(\mathtt{A})}))$ with the property
$t_{j-1}\text{Cone}(\psi_{j})[-1] =
t_{j-1}\mathcal{Y}_{j-1}^{\bullet}$. Set $\mathcal{Y}_{j}^{\bullet} 
= \text{Cone}(\psi_{j})[-1]$ and the induction follows.  

Define the \textit{limit complex} $\mathcal{Y}^{\bullet} \in
\Ka(\add(Y^{(\mathtt{A})}))$ by $t_j\mathcal{Y}^{\bullet} =
t_j\mathcal{Y}^{\bullet}_j$ for all $j \geq 0$ (that  
for all $ k < j$ we have $t_k\mathcal{Y}_{j}^{\bullet} = 
t_k\mathcal{Y}_{k}^{\bullet}$, guarantees that 
$\mathcal{Y}^{\bullet}$ is well-defined). Moreover, we have a 
quasi-isomorphism $\Phi:\mathcal{Y}^{\bullet} \to 
\mathcal{X}^{\bullet}$. Hence the general statement is true and we see
by the construction that  $\mathcal{Y}^{\bullet}$ is bounded,
whenever $\mathcal{X}^{\bullet}$ is.
\end{proof}

\begin{lemma}\label{apl2}
Let $\mathtt{A}$ be a finite-dimensional associative algebra, $X^{(\mathtt{A})}$ 
be a $\mathtt{A}$-module and $Y^{(\mathtt{A})}$ be a (generalized) (co)tilting 
$\mathtt{A}$-module.  Assume that $X^{(\mathtt{A})}$ admits a (possibly infinite) 
resolution by modules from
$\mathrm{add}(Y^{(\mathtt{A})})$. Then for every negative complex
$\mathcal{X}^{\bullet}\in \Ka^b(\add(X^{(\mathtt{A})}))$ there exists
a (possibly infinite) negative complex, $\mathcal{Y}^{\bullet}\in 
\Ka(\add(Y^{(\mathtt{A})}))$, which is quasi-isomorphic to 
$\mathcal{X}^{\bullet}$.
\end{lemma}

\begin{proof}
The statement is proved by induction analogous to that used in the previous 
lemma. Moreover, since we start with a finite complex from the very 
beginning, no truncation is needed. We leave the details out.
\end{proof}

%%%%%%%%%%%%%%%%%%%%%%%%%%%%%%%%%%%%%%%%%%%%%%%%%%%%%%%%%%%%%%%%%%%%%%%%%%%
\vspace{1cm}
%%%%%%%%%%%%%%%%%%%%%%%%%%%%%%%%%%%%%%%%%%%%%%%%%%%%%%%%%%%%%%%%%%%%%%%%%%%

\begin{flushleft}
\bf Acknowledgments. 
\end{flushleft}

The research was partially supported by The Swedish Foundation for
International Cooperation in Research and Higher Education. 
For the second author the research was also partially supported by The
Royal Swedish Academy of Science and by The Swedish Research Council. 
We would like to thank the referee for a very careful reading of the 
paper and for many remarks and suggestions which led to the improvements
in the paper.

\vspace{1cm}

\noindent
Department of Mathematics, Uppsala University, Box 480, SE-75106,
Uppsala, SWEDEN, 

\noindent
{\tt frisken\symbol{64}math.uu.se}\\
{\tt mazor\symbol{64}math.uu.se}

\end{document}